\documentclass[twocolumn]{svjour3}          % twocolumn
\smartqed  % flush right qed marks, e.g. at end of proof
\usepackage{graphicx}
\usepackage{amssymb,amsmath}
\usepackage{graphicx}
\usepackage{color}
\usepackage{multirow}
\usepackage{todonotes}

\newcommand{\eps}{\varepsilon}
\newcommand{\R}{\mathbb{R}}

\renewcommand{\vec}[1]{\mathbf{#1}}

\journalname{Biological Cybernetics}

\begin{document}

\title{Excitable Networks for Finite State Computation with Continuous Time Recurrent Neural Networks}

\author{Peter Ashwin         \and
        Claire Postlethwaite %etc.
}

\institute{Peter Ashwin \at
              Center for Systems, Dynamics and Control,  Department of Mathematics, University of Exeter, Exeter EX4 4QF, UK
              \email{p.ashwin@exeter.ac.nz}           %  \\
           \and
          Claire Postlethwaite \at
          Department of Mathematics, University of Auckland, Auckland, 1142, New Zealand.
}

\date{Received: date / Accepted: date}

\maketitle

\begin{abstract}
Continuous time recurrent neural networks (CTRNN) are systems of coupled ordinary differential equations that are simple enough to be insightful for describing learning and computation, from both biological and machine learning viewpoints. We describe a direct constructive method of realising finite state input-dependent computations on an arbitrary directed graph. The constructed system has an excitable network attractor whose dynamics we illustrate with a number of examples. The resulting CTRNN has intermittent dynamics: trajectories spend long periods of time close to steady-state, with rapid transitions between states. Depending on parameters, transitions between states can either be \emph{excitable} (inputs or noise needs to exceed a threshold to induce the transition), or \emph{spontaneous} (transitions occur without input or noise).  In the excitable case, we show the threshold for excitability can be made arbitrarily sensitive.
\keywords{
Continuous Time Recurrent Neural Network \and Nonlinear Dynamics \and Excitable Network Attractor }
\end{abstract}

\section{Introduction}
\label{sec:intro}

It is natural to try to understand computational properties of neural systems through the paradigm of network dynamical systems, where a number of dynamically simple units (i.e. with attracting equilibria or periodic orbits) interact to give computation as an emergent property of the system. This is the as much the case for biological models of information processing, pattern generation and decision making as it is for artificial neural networks inspired by these. A variety of specific models have been developed to describe the dynamics and training of recurrent networks comprised of coupled neurons in biological and artificial settings. The particular challenge that we address here is the construction of arbitrarily complex, but specified, dynamical structures that enable discrete (finite-state) computation in an input-driven system that is nonetheless continuous in both state and time. Clearly, invariant objects of the autonomous system (such as equilibria, periodic orbits and chaotic attractors) only form part of the picture and an input-dependent (non-autonomous) approach such as \cite{manjunath2012theory} is needed. 

To help understand the response of systems to inputs the authors introduced in~\cite{Ashwin2016a} a notion of a ``network attractor'', namely an invariant object in phase space that may contain several local invariant sets, but also systems of interconnections between them. This generalises the notion of ``heteroclinic network'' and ``winnerless competition/stable heteroclinic channels'' \cite{afraimovich2004origin} which have been used to describe a range of sequence generation, computational and cognitive effects in neural systems: see for example \cite{rabinovich2001dynamical,rabinovich2006generation,afraimovich2004heteroclinic,rabinovich2020sequential,hutt2017sequences}. These models connect computational states (represented as saddles in the dynamics). In the presence of inputs or noise, the switching between saddles is a useful model for spontaneous computational processes but there are two problems. One of these is that the states are dynamically unstable and so there are spontaneous transitions that are not noise driven. The other is that heteroclinic chains are destroyed by arbitrarily small perturbations, unless there are special structures in phase space (symmetries or invariant subspaces).

Network attractors joining stable equilibria \cite{Ashwin2016a} overcome both of these issues. In this paper we show that arbitrarily complex network attractors exist in idealised class of model is the continuous time recurrent neural network (CTRNN) \cite{beer1995dynamics}. A CTRNN is a set of differential equations each with one scalar variable that represents the level of activation of a neuron, and feedback via a saturating nonlinearity or ``activation function''. These models have been extensively investigated in the past decades as simple neurally-inspired systems that can nonetheless (without input) have complex dynamics by virtue of the nonlinearities present \cite{beer1995dynamics}. They can also be trained to perform arbitrarily complex tasks depending on time-varying input \cite{tuci2002evolutionary,yamauchi1994sequential}. They are frequently used in investigations of evolutionary robotics~\cite{blynel2003exploring}  and (in various equivalent formulations \cite{chow2019before}) as models for neural behaviour in biological or cognitive systems. For example, the classical work of Hopfield and Tank \cite{hopfield1985neural} considers such systems with symmetric weights to solve optimization problems, while more recently,~\cite{bhowmik2016reservoir} discusses CTRNN models for episodic memory, and several other biological and cognitive applications are discussed in \cite{nikiforou2019dynamics}.
 
CTRNNs are often referred to as ``universal dynamical approximators'' \cite{funahashi1993approximation}, meaning that the trajectory of a CTRNN can approximate, to arbitrary precision, any other prescribed (smooth) trajectory in $\R^n$.  However, this does not mean that the dynamics of CTRNNs are simple to understand, or that it is easy to form the above approximation. It also raises the question of how a ``trained" CTRNN performs a complex task. Gradient descent or more general evolutionary training algorithms train the network by navigating through a high dimensional landscape of possible feedback weightings and moving these towards a setting that is sufficiently optimal for the task. It may be possible to give a clear description of the resulting nonlinear dynamics of the autonomous (constant input) CTRNN, but we want to understand not only this but also how inputs affect the state of the system. 
 
The main theoretical result in this paper focusses on ``excitable network attractors'': these consist of a finite set of local attracting equilibria and excitable connections between them (see Appendix \ref{app:networks}). It was demonstrated in \cite{Ashwin2018sensitive} that an excitable network attractor can be used to embed an arbitrary Turing machine within a class of purpose-designed coupled dynamical system with two different cell types. Rather than relying on an optimization approach to design the system, that paper gave a \emph{constructive} method for designing a realisation of any desired network attractor. However, this construction required specialist dynamical cells with quite complex nonlinear couplings between them, and a comparatively large number of cells. It was recently shown in \cite{ceni2020interpreting} that trained RNNs in the form of echo-state networks can realise excitable networks for certain tasks and that structural errors in the trained network can explain errors in imperfectly trained systems.

The current paper demonstrates that CTRNN dynamics is sufficiently rich to realise excitable networks with arbitrary graph topology, simply by specifying appropriate connection weights. The construction algorithm in the proof of Theorem~\ref{thm:main}, assigns one of only four values to each of the connection weights to realise an arbitrary graph on $N$ vertices as an excitable network on $N$ states (subject to some minor constraints on its connectivity), using a CTRNN with $N$ cells. The CTRNN we consider in this paper (see for example \cite{beer1995dynamics}, which corresponds to a continuous time Hopfield model \cite{hopfield1985neural} in the symmetric coupling case $w_{ij}=w_{ji}$) are ordinary differential equations 
\begin{equation}
\dot{y}_i = -y_i+\sum_{j} w_{ij} \phi(y_j)+I_i(t),
\label{eq:ctrnn_ode}
\end{equation}
where $\vec{y}=(y_1,\ldots y_N)\in \R^N$ is the internal state of the $N$ cells of the system, $w_{ij}$ is a matrix of connection weights, $\phi$ is a (sigmoid) activation function that introduces a saturating nonlinearity into the system and $I_i(t)$ is an input. We say system (\ref{eq:ctrnn_ode}) is {\em input-free} if $I_i(t)=0$ for all $i$ and $t$.

We consider two cases for $\phi$, a smooth function
\begin{equation}
\label{eq:phi}
\phi(y)=\phi_S(y):=\left[1+\exp\left(-\frac{(y-\theta)}{\epsilon}\right)\right]^{-1},
\end{equation}
and a piecewise affine function
\begin{equation}
\label{eq:phip}
\phi(y)=\phi_P(y):=\begin{cases}
0 & y-\theta <-2\epsilon,\\
(y-\theta)/(4\epsilon)+1/2 & |y-\theta|\leq 2\epsilon,\\
1 & y-\theta>2\epsilon,
\end{cases}
\end{equation}
In both cases, $\phi$ is monotonic increasing with
$$
\lim_{y\rightarrow\infty}\phi(y)=1, ~ \lim_{y\rightarrow-\infty}\phi(y)=0, ~ 
\phi(\theta)=1/2,
$$
and a maximum derivative at $y=\theta$, equal to 
$
\frac{1}{4\epsilon}.
$
In both cases, $\epsilon$ and $\theta$ are parameters, and we are interested in the case $0<\epsilon \ll 1$. 
In general, the function $\phi$ need not be the same in every component of~\eqref{eq:ctrnn_ode}, but here we make a simplifying assumption that it is.

Note that both activation functions (\ref{eq:phi}) and (\ref{eq:phip}) have piecewise constant limits in the singular limit $\epsilon\rightarrow 0$. Such limiting systems are of Fillipov type and have been explored in various biological contexts, especially for gene regulatory dynamics. These can also have rich dynamics as discussed in the literature (for example \cite{gouze2003aclass,harris2015bifurcations}) but we do not consider this limit here.

The main contribution of Section~\ref{sec:method} is to give, in  Theorem~\ref{thm:main}, a construction of a connection weight matrix $w_{ij}$ such that the dynamics of the input-free system~\eqref{eq:ctrnn_ode} contains (or \emph{realises}: definition given below) an excitable network attractor, as defined in~\cite{Ashwin2016a}. We prove this (with details in Appendix~\ref{app:proof}) for the case of the piecewise affine function $\phi_P$. In section~\ref{sec:ex} we present evidence that this is also true for the smooth case $\phi_S$ for an open set of parameters. Qualitatively, this means the system will contain a number of stable equilibrium states, and small inputs (either deterministic, or noisy) will push the trajectory from one stable equilibrium into the basin of attraction of another. In this way, transitions can be made around the network, and the transition time between states tends to be much smaller than the residence times of the trajectory in neighbourhoods of the states. In particular, we can choose $w_{ij}$ so that the network attractor has (almost) any desired topology. In Appendix~\ref{app:networks} we recall formal definitions of network attractors from~\cite{Ashwin2016a,Ashwin2018sensitive}.

In section~\ref{sec:ex} we consider several examples of simple graphs, and demonstrate that the desired networks do indeed exist in the systems as designed. We also perform numerical bifurcation analysis to demonstrate the connection between periodic orbits in the input-free deterministic system~\eqref{eq:ctrnn_ode} and excitable networks in the same system with additive noise, that is, the system of stochastic differential equations (SDEs):
\begin{equation}
dy_i = \left(-y_i+\sum_{j} w_{ij} \phi(y_j) \right)dt+ \sigma dW_i(t),
\label{eq:ctrnn_sde}
\end{equation}
where $W_i(t)$ are independent standard Wiener processes. Here, the noise plays the role of inputs that propel the trajectory around the network, although of course this occurs in a random manner. In sections~\ref{sec:KS} and~\ref{sec:tennode} we consider graphs that have multiple edges leading out from a single vertex, and show that additional equilibria may appear in the network attractor where two or more cells are active simultaneously. We further show that the existence of these additional equilibria can be suppressed by choosing one of the parameters used in the construction of the weight matrix $w_{ij}$ to be sufficiently large.

Section~\ref{sec:discuss} concludes by relating our results to other notions of sequential computation. We also conjecture some extensions of the results shown in this paper.

\section{Construction of a CTRNN with a network attractor}
\label{sec:method}

Let  $G$ be an arbitrary directed graph  between $N$ vertices, and let  $a_{ij}$ be the adjacency matrix of $G$. That is, $a_{ij}=1$ if there is a directed edge from vertex $i$ to vertex $j$, and $a_{ij}=0$ otherwise. 

Let $\Sigma$ be an invariant set for a system of ordinary differential equations. 
We say $\Sigma$ is an excitable network that \emph{realises} a graph $G$ for some amplitude $\delta>0$ if for each vertex $v_i$ in $G$ there is a unique stable equilibrium $\xi_i$ in $\Sigma$, and if there is an excitable connection with amplitude $\delta$ in $\Sigma$ from $\xi_i$ to $\xi_j$ whenever there is an edge in $G$ from $v_i$ to $v_j$. The existence of an excitable connection means that there exists a trajectory with initial condition within a distance $\delta$ of $\xi_i$, which asymptotes in forward time to $\xi_j$
%We say the realisation is \emph{proper} if there is a one-to-one correspondence between edges in $G$ and excitable connections 
(formal definitions are given in Appendix~\ref{app:networks}).

\begin{figure}
	\begin{center}
		\includegraphics[width=0.45\textwidth]{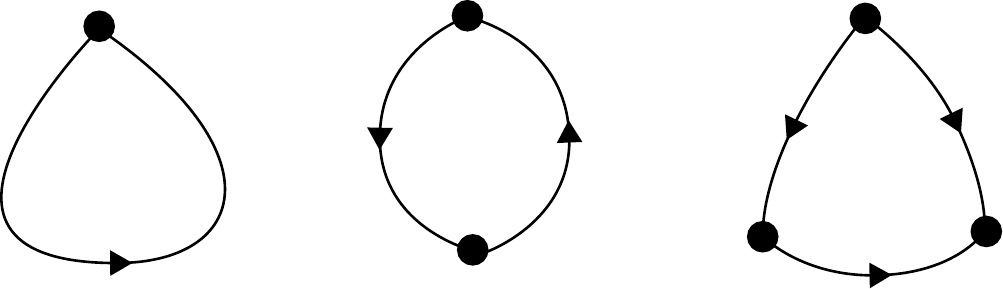}
	\end{center}
	\caption{\label{fig:graph_constraints} From left to right, the figures show an order-one loop, an order-two loop, and a $\Delta$-clique in a directed graph. The construction we provide realises an arbitrary graph $G$, as long of none of these subgraphs are present in $G$.}
\end{figure}

For the purposes of our construction of a network attractor, we assume that $G$ contains no loops of order one, no loops of order two, and no $\Delta$-cliques. Figure~\ref{fig:graph_constraints} shows each of these graph components schematically. 
In terms of the adjacency matrix, for $G$ to contain no loops of order one requires that 
\begin{equation}
\label{eq:aijassumptions1}
a_{ii}=0 \mbox{ for all } i;
\end{equation}
for $G$ to contain no loops of order two requires that
\begin{equation}
\label{eq:aijassumptions2}
a_{ij}a_{ji}=0 \mbox{ for all }i,j;
\end{equation}
and  for $G$ to contain no $\Delta$-cliques requires that
\begin{equation}
\label{eq:aijassumptions2}
a_{ij}a_{ik}a_{jk}=0 \mbox{ for all }i,j,k.
\end{equation}
In our earlier work, we have demonstrated a network design which can admit order-two loops and $\Delta$-cliques~\cite{Ashwin2016a}, although this previous construction is not motivated by neural networks \emph{per se}, and requires a higher dimensional system of ODEs (for a given graph) than the one presented here.

Before we move into the details of the construction, we briefly discuss our terminology. A graph $G$ has \emph{vertices}, which correspond to (stable) \emph{equilibria} in the phase space of the dynamical system~\eqref{eq:ctrnn_ode}. Also within the phase space, there exist \emph{excitable connections} (sometimes abbreviated to \emph{connections}) between the equilibria, which correspond to the \emph{edges} of the graph. When a trajectory in the phase space moves between neighbourhoods of equilibria along a path close to one of these connections, we say that a \emph{transition} between the equilibria has occurred. We refer to each of the components of the dynamical system~\eqref{eq:ctrnn_ode} as a \emph{cell}, and say that a cell $j$ is \emph{active} if $\phi(y_j)$ is close to one.

\subsection{Realization of arbitrary directed graphs as network attractors}
\label{sec:thm}

We construct  a weight matrix $w_{ij}$ that depends only on the adjacency matrix $a_{ij}$, and on four parameters $w_t$, $w_s$, $w_m$ and $w_p$. It is given by
\begin{equation}\label{eq:w}
w_{ij}=w_t+(w_s-w_t)\delta_{ij}+(w_p-w_t)a_{ji}+(w_m-w_t)a_{ij},
\end{equation}
where $\delta_{ij}$ is the Kronecker $\delta$.
This choice of $w_{ij}$ ensures that $w_{ii}=w_s$, $w_{ij}=w_p$ if there is a directed edge in $G$ from vertex $j$ to $i$ (i.e. $a_{ji}=1$), $w_{ij}=w_m$ if there is a directed edge in $G$ from vertex $i$ to vertex $j$ (i.e. $a_{ij}=1$), and $w_{ij}=w_t$ otherwise. In later sections we allow for different weights along different edges by allowing $w_p$ to depend on $i$ and $j$ (i.e. $w_p=w_p^{ij}$). We give an overview of how each of the parameters affect the dynamics of the system in section~\ref{sec:smooth}.

We write $\vec{w}=(\epsilon,\theta,w_s,w_m,w_t,w_p)\in \R^{6}$ to be a vector of all parameter values. The next result shows that for the piecewise affine activation function and suitable choice of parameters, there is an embedding of $G$ as an excitable network attractor for the input-free system.

\begin{theorem}\label{thm:main}
	For any directed graph $G$ with $N$ vertices containing no loops of order one, no loops of order two, and no $\Delta$-cliques, and any small enough $\delta>0$, there is an open set $W_{\mathrm{ex}}\subset\R^{6}$ such that if the parameters $\vec{w}\in W_{\mathrm{ex}}$ then the dynamics of input-free equation~\eqref{eq:ctrnn_ode} on $N$ cells with piecewise affine activation \eqref{eq:phip} and $w_{ij}$ defined by~\eqref{eq:w} contains an excitable network attractor with threshold $\delta$, that realises the graph $G$.
\end{theorem}

Recall that by \emph{realises} we mean that all edges in the graph are present as transitions between stable equilibria using perturbations of size at most $\delta$.

\smallskip

\proof
We give the main ideas behind the proof here, deferring some of the details to Appendix \ref{app:proof}.  We construct an excitable network attractor in $\R^{N}$ for \eqref{eq:ctrnn_ode} with piecewise activation function \eqref{eq:phip} and weight matrix \eqref{eq:w}. For any $\frac{1}{2}>\delta>0$, we show there exist parameters $\bf{w}$ (with $\epsilon>0$ small) and stable equilibria $\xi_k$ ($k=1,\ldots,N$) that are connected according to the adjacency matrix $a_{ij}$ by excitable connections with amplitude $\delta$. Below, we provide an explicit set of parameters that make such a realisation, and note that the realisation will hold for an open set of nearby parameters.

We show in Appendix \ref{app:proof} that the equilibria $\xi_k$ have components (cells) that are close to one of four values $Y_T$, $Y_D$, $Y_L$, $Y_A$ related to the edges attached to the corresponding vertex $k$ in the graph $G$. For any $0<\delta<\frac{1}{2}$, we use the following parameters:
\begin{equation}
\label{eq:paramchoicea}
\epsilon=\dfrac{\delta}{8},~\theta=\frac{1}{2},~w_s=1,~w_t=0, 
\end{equation}
and then $w_p$ and $w_m$ are given by
\begin{equation}
\label{eq:paramchoiceb}
w_p=\theta-\dfrac{\delta}{2},\quad w_m=-(w_s-\theta)-\dfrac{\delta}{2}.
\end{equation}
We then set
\begin{equation}\label{eq:yaltd}
\begin{split}
Y_A:=w_s,\ Y_L:=w_p=\theta-\frac{\delta}{2},\  \\
Y_T:=w_m=-(w_s-\theta)-\frac{\delta}{2},\  Y_D:=w_t
\end{split}
\end{equation}

We use square brackets and subscripts to identify the components of points in phase space, that is, $[\xi_k]_j$ is the $j$th component of $\xi_k$. Each $\xi_k$ has:
\begin{itemize}
	\item Exactly one cell that is {\em Active}: $[\xi_k]_k= Y_A$
	\item A number of cells that are {\em Leading}: $[\xi_k]_j= Y_L$ (if $a_{kj}=1$)
	\item A number of cells that are {\em Trailing}: $[\xi_k]_j= Y_T$ (if $a_{jk}=1$)
	\item All remaining cells are {\em Disconnected}: $[\xi_k]_j= Y_D$ ($a_{kj}=a_{jk}=0$).
\end{itemize}

Note that the requirement of no loops of order one or two and no $\Delta$-cliques implies that this labelling is well defined. 

\begin{figure*}
	\centering
	\includegraphics[width=13cm]{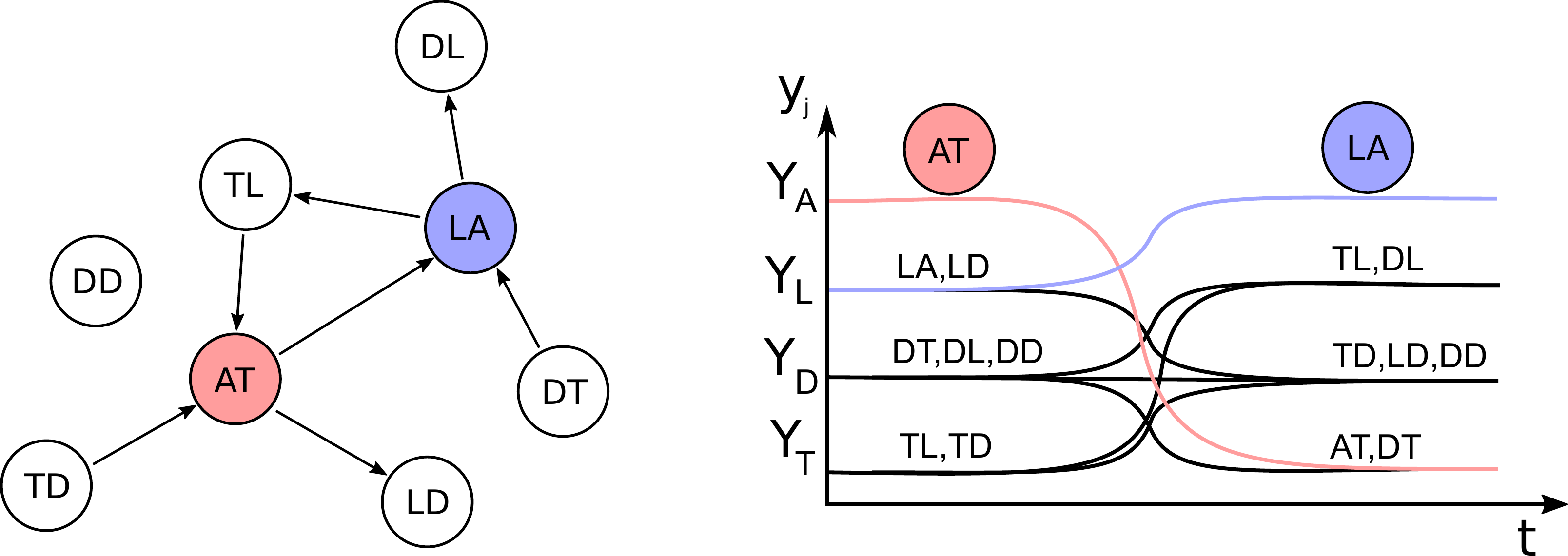}
	\caption{Schematic diagram showing a transition in the network attractor of Theorem~\ref{thm:main}. Left: a schematic of the directed graph showing the eight distinct types of vertex. Right: schematic time series; as the Active cell changes there will be a transition in the connected Leading and Trailing cells as shown.  Which cell becomes active is controlled by a small perturbation to the set of currently Leading cells. 
	}
	\label{fig:transition_schematic}
\end{figure*}

From  equilibrium $\xi_k$, there is an excitable connection to any of the equilibria $\xi_l$ with $a_{kl}=1$, that is, any of the Leading cells can become the Active cell. During a transition, the remaining cells can be classified into six types, which are identified in figure~\ref{fig:transition_schematic}, and depend (for each $j$) on the values of the four entries in the adjacency matrix  $a_{jk}$, $a_{kj}$, $a_{jl}$ and $a_{lj}$. We label cell $k$ as AT (Active--Trailing) and cell $l$ as LA (Leading--Active).
Note that the lack of two cycles means that the cases with $a_{jk}=a_{kj}=1$ or $a_{jl}=a_{lj}=1$ (a total of seven possibilities) cannot occur, and the lack of $\Delta$-cliques mean that the cases with $a_{jk}=a_{jl}=1$, $a_{kj}=a_{jk}=1$, or $a_{kj}=a_{lj}=1$ also cannot occur (which includes the cases where a cell would switch from Leading to Trailing). The remaining six possibilities are listed below.

\begin{itemize}
\item Type DD: $a_{jk}=a_{kj}=a_{jl}=a_{lj}=0$; the cell is Disconnected throughout.
\item  Type TD: $a_{jk}=1$, $a_{kj}=a_{jl}=a_{lj}=0$; the cell switches from Trailing to Disconnected.
\item Type LD: $a_{kj}=1$, $a_{jk}=a_{jl}=a_{lj}=0$; the cell switches from Leading to Disconnected.
\item Type TL: $a_{jk}=a_{lj}=1$, $a_{kj}=a_{jl}=0$; the cell switches from Trailing to Leading.
\item Type DT: $a_{jl}=1$, $a_{jk}=a_{kj}=a_{lj}=0$; the cell switches from Disconnected to Trailing.
\item Type DL: $a_{lj}=1$, $a_{jk}=a_{kj}=a_{jl}=0$; the cell switches from Disconnected to Leading.
\end{itemize}
The right panel in figure~\ref{fig:transition_schematic} shows how a transition from cell AT active to cell LA active will occur in a general network. 

To prove the existence of an excitable connection giving a realisation we consider a perturbation from $\xi_k$ to the point
$$
\zeta_{k,l}=\xi_k+\delta e_l,
$$
where $e_l$ is the unit basis vector, and we show in Appendix~\ref{app:proof} that, for small enough $\delta$, $\zeta_{k,l}$ is in the basin of attraction of $\xi_l$ if $a_{kl}=1$. This means there is an excitable connection from $\xi_k$ to $\xi_l$ in this case.
\qed

~

We believe that for small enough $\delta$ and suitable choice of weights, the realisation of $G$ in Theorem~\ref{thm:main} can be made {\em almost complete} in the sense analogous to the similar notion for heteroclinic networks \cite{ashwin2020almost}, namely that the set
$$
\bigcup_{\xi\in E} B_{\delta}(\xi) \setminus \Sigma_E
$$
has zero measure. If a network realization is almost complete then almost all trajectories starting close to some $\xi_k$ will remain close, or will follow a connection corresponding to the realization.  We do not have a proof of this, though Appendix~\ref{sec:proofproper} shows that for small enough $\delta$, $\zeta_{k,l}$ is in the basin of attraction of $\xi_l$ if and only if $a_{kl}=1$. This does not preclude the possibility that there exist perturbations in $B_{\delta}(\xi_i)$ other than $\zeta_{k,l}$ resulting in trajectories asymptotic to $\xi_l$, or indeed to other attractors. Our numerical investigations suggest that by choosing $w_t$ non-zero and large enough, any connections to other equilibria can be suppressed. This is discussed more in sections~\ref{sec:tennode} and~\ref{sec:discuss}.

\subsection{Excitable networks for smooth nonlinearities}
\label{sec:smooth}

For small enough $\epsilon$, the smooth activation function \eqref{eq:phi} can be made arbitrarily close to the piecewise activation function~\eqref{eq:phip}, and so we expect Theorem~\ref{thm:main} to also apply in the smooth case, but do not give a proof here. However, throughout section~\ref{sec:ex} we use the smooth activation function~\eqref{eq:phi} in our examples. We use
\begin{equation}\label{eq:defparms}
\begin{split}
&\epsilon=0.05,~\theta=0.5,~w_s=1, \\ &w_m=-0.7,  w_p=0.3,~w_t=0
\end{split}
\end{equation}
as our default parameter set (compare this choice of parameters with that given in equations~\eqref{eq:paramchoicea} and~\eqref{eq:paramchoiceb}, with $\delta=0.4$), though there will be an open set of nearby parameters with analogous behaviour. In section~\ref{sec:ex} we provide several examples of using this choice of weight matrix to realise a graph $G$.

From equation~\eqref{eq:yaltd} we can see that the parameter choices directly affect the location of the equilibria $\xi_k$ in phase space. As we will see in the following sections, the parameters also have further effects on the dynamics. In particular, the relative sizes of the parameters $w_p$ and $\theta$ determine whether the dynamics are excitable or spontaneous: essentially, for $\epsilon$ small enough, $w_p$ needs to be smaller than $\theta$ to observe excitable dynamics. If $w_p$ is too large, then the equilibria $\xi_k$ cease to exist: periodic orbits can exist instead. The parameter $w_m$ controls how fast a Trailing cell decays, and the parameter $w_t$ controls the suppression effects when there is more than one Leading cell. We discuss the effects of $w_t$ in more detail in section~\ref{sec:KS}.

For a smooth activation function such as (\ref{eq:phi}) that is invertible on its range, there is a useful change of coordinates to $J_i=\phi(y_i)$ (a similar transformation is made in~\cite{beer1995dynamics}). The input-free equations~\eqref{eq:ctrnn_ode} then become
\begin{equation}
\dot{J}_i = \frac{1}{\epsilon}J_i(1-J_i)\left(-\phi^{-1}(J_i)+\sum_{j} w_{ij} J_j\right),
\label{eq:ctrnn_Jode}
\end{equation}
where each $J_i\in (0,1)$ (which is the domain of the function $\phi^{-1}$), and
\[
\phi^{-1}(x)=\theta-\epsilon \ln\left(\frac{1-x}{x}\right).
\]

Each vertex in the graph $G$ is realised in the phase space of the $J_k$ variables by a stable equilibrium with one of the $J_k$ close to $1$ and the remainder close to $0$. With a slight abuse of notation, we still refer to these equilibria as $\xi_k$.
As we will see in the examples which follow, the parameters are chosen such that the dynamics are close to a saddle-node bifurcation; more precisely such that the system is near a codimension $N$ bifurcation where there are $N$ saddle-nodes of the equilibria $\xi_k$. This bifurcation has a global aspect analogous to a saddle-node on an invariant circle (SNIC) or saddle-node homoclinic bifurcation in that coupling weights are such that there are global connections between the saddle nodes that reflect the network structure.

If parameters are on one side of the $N$ saddle-node bifurcations then for each $k$ there will exist a pair of equilibria, one of which will be stable (this is the equilibrium $\xi_k$), and a nearby saddle.  A small perturbation can then move the trajectory out of the basin of attraction of the stable equilibrium and effect a transition to another equilibrium: this is an \emph{excitable} connection from $\xi_k$. 

If parameters are on the other side of the $N$ saddle-node bifurcations then the flow through the corresponding $N$ region of phase space will be slow but we will still observe intermittent dynamics as the trajectory passes through this region \cite[p99]{strogatz1994nonlinear}.
In this case, we refer to a region where (a) there is a unique local minimum of $|\dot{\vec{y}}|$ and (b) a large subset of initial conditions in this region pass close to this minimum, as a \emph{bottleneck region} $P_k$, and refer to a \emph{spontaneous} transition past $P_k$. (This is also called the \emph{ghost of a saddle-node bifurcation} in \cite{strogatz1994nonlinear}.) If all the connections corresponding to edges in the graph $G$ are excitable, then the system contains an excitable network attractor. If all connections are spontaneous, then we typically see a periodic orbit, although we do not prove this.

\section{Examples for the smooth activation function}
\label{sec:ex}

\subsection{Two vertex graph}
\label{sec:twonode}

For our first example we consider the connected graph with two vertices and a single edge joining them, that is $a_{12}=1$, and $a_{ij}=0$ for $(i,j)\neq (1,2)$. We use bifurcation analysis to show that the transition between spontaneous and excitable dynamics is caused by a saddle-node bifurcation, and find an approximation to the location of the saddle-node bifurcation in parameter space.

\begin{figure*}
	\setlength{\unitlength}{0.8mm}
	\begin{center}
		\begin{picture}(160,75)(0,0)
		\put(5,5){\includegraphics[width=5.6cm]{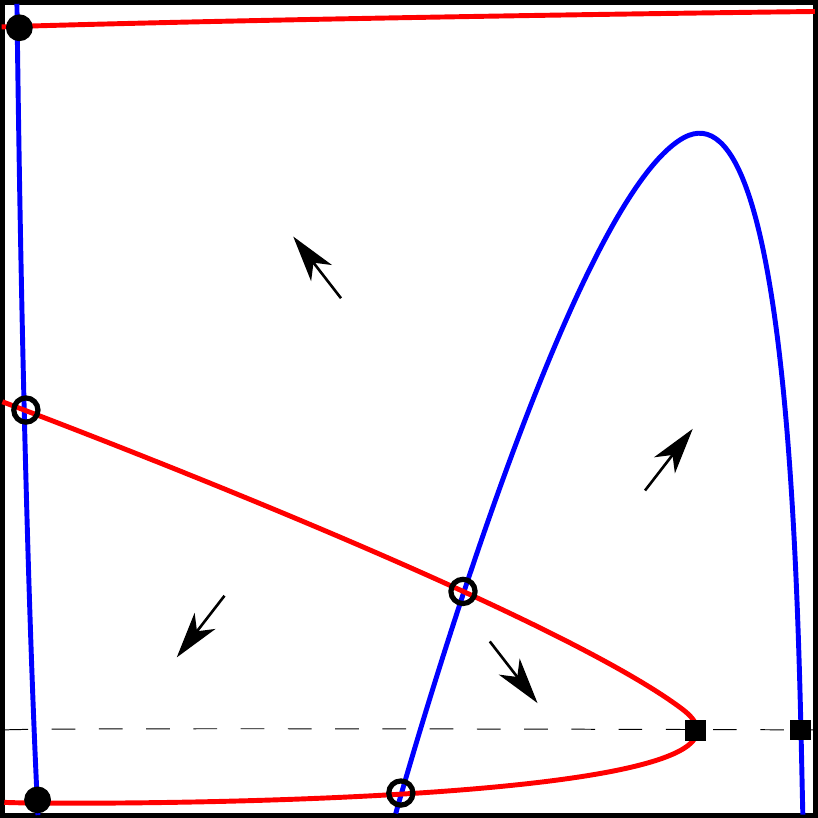}} %5.6=7*0.8
		\put(85,5){\includegraphics[width=5.6cm]{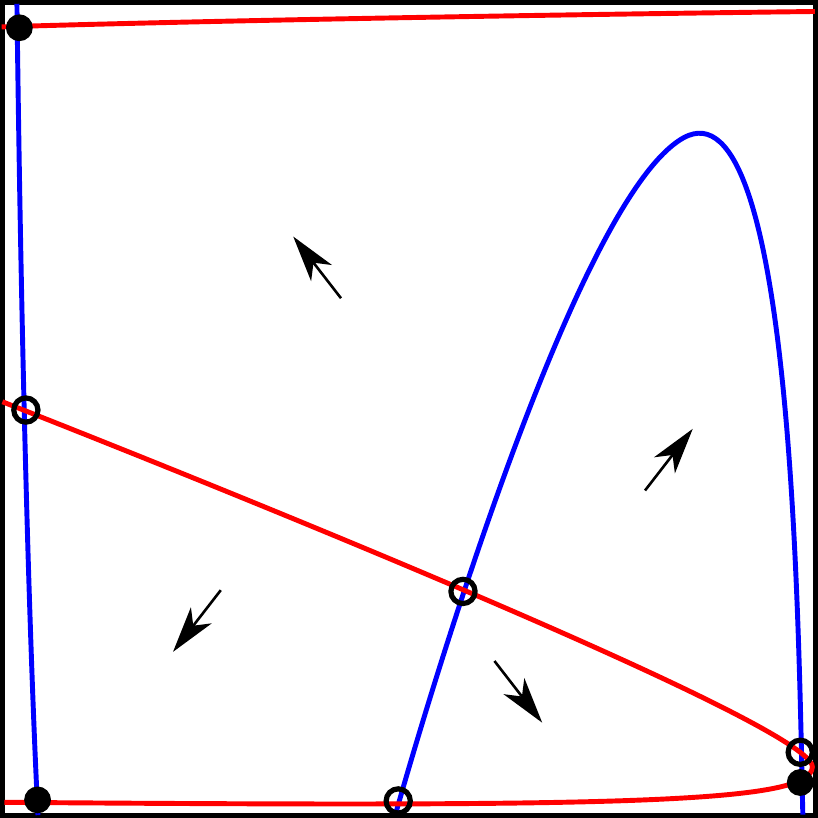}}
		\put(69,1){$J_1$}
		\put(149,1){$J_1$}
		\put(-1,69){$J_2$}
		\put(79,69){$J_2$}
		\put(-1,12){$x_-$}
		\put(64,15){$A$}
		\put(69,15){$B$}
		\end{picture}
	\end{center}
	\caption{\label{fig:jnull} Schematic illustration of the phase space and nullclines of the system~\eqref{eq:J2odes}.  The blue curves are the $J_1$ nullclines, and the red curves are the $J_2$ nullclines. Solid dots show stable equilibria and open dots show unstable equilibria. The parameter $w_p$ is decreased from left figure to right figure, creating two further equilibrium in a saddle-node bifurcation. The squares labelled $A$ and $B$ are referred to in the proof of lemma~\ref{lem:wpsn}.  }
\end{figure*}
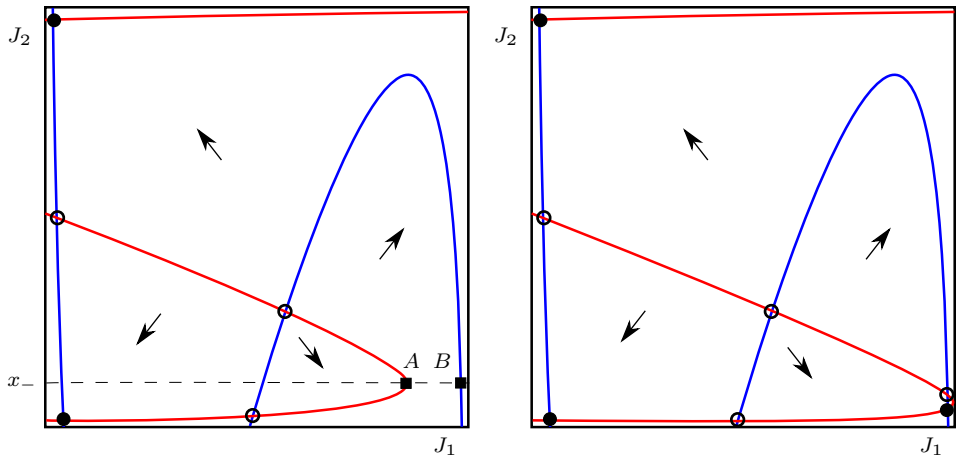

The two-dimensional system of equations is:
\begin{align}
\dot{y}_1&= -y_1+w_s\phi(y_1)+w_m\phi(y_2),\\
\dot{y}_2&= -y_2+w_s\phi(y_2)+w_p\phi(y_1). 
\end{align}
In the $J_i$ variables, this becomes
\begin{align} \label{eq:J2odes}
\dot{J}_1&=\frac{1}{\epsilon}J_1(1-J_1)\left(-\phi^{-1}(J_1)+w_sJ_1+w_mJ_2\right), \\
\dot{J}_2&=\frac{1}{\epsilon}J_2(1-J_2)\left(-\phi^{-1}(J_2)+w_sJ_2+w_pJ_1\right), \nonumber
\end{align}
where $(J_1,J_2)\in (0,1)^2$.
We note the following properties of the function $g:(0,1)\rightarrow\R$, with $g(x)=\phi^{-1}(x)-w_s x$:
\[
g'(x)=\frac{\epsilon}{x(1-x)}-w_s,\quad g''(x)=\frac{\epsilon(2x-1)}{x^2(1-x)^2},
\]
\[\lim_{x\rightarrow 0}g(x)=-\infty,\quad  \lim_{x\rightarrow 1}g(x)=\infty,\]
\[  g\left(\frac12\right)=\theta-\frac{w_s}{2},\quad g'\left(\frac12\right)=4\epsilon-w_s.\] 
If $w_s>4\epsilon$, then $g$ has local extrema at $x_+$ and  $x_-$, where
\begin{equation}\label{eq:xpm}
x_{\pm}=\frac{1}{2}\pm\sqrt{\frac{1}{4}-\frac{\epsilon}{w_s}}.
\end{equation}
The $J_1$ and $J_2$ nullclines of system~\eqref{eq:J2odes} are at, respectively
\begin{equation}\label{eq:nullclines}
J_2=\frac{1}{w_m}g(J_1) \quad\text{and} \quad J_1=\frac{1}{w_p}g(J_2).
\end{equation}
In figure~\ref{fig:jnull}, we show a sketch of the phase space of~\eqref{eq:J2odes}; the $J_1$ nullcline is shown in blue, and the $J_2$ nullcline in red. Solid dots show stable equilibria, open dots show unstable equilibria and arrows show the direction of flow. Note that we do not include any nullclines at $J_i=0$ or $J_i=1$ because they are not in the domain of equations~\eqref{eq:J2odes}. As $w_p$ is decreased (in the figures, moving from left to right), a saddle-node bifurcation creates a pair of equilibrium solutions. 

\begin{lemma} \label{lem:wpsn}
	If $\theta<w_s$, and $0<\epsilon \ll 1$, then a saddle node bifurcation occurs in system~\eqref{eq:J2odes} when 
	\begin{equation}\label{eq:wpsn}
	w_p=w_p^{SN}\equiv\epsilon \log\epsilon +\theta-\epsilon (1+\log w_s)+\frac{\epsilon^2}{w_s}+O(\epsilon^3).
	\end{equation}
\end{lemma}

To begin the proof, we note that the saddle-node bifurcation will occur when the points $A$ and $B$ (marked by squares in the left-hand panel of figure~\ref{fig:jnull}) coincide. These points are defined as the intersection of the nullclines with the line at $J_2=x_-=\frac{\epsilon}{w_s}+O(\epsilon^2)$, i.e. at the local extrema of the $J_2$ nullcline. 

Let the $J_1$ coordinate of $A$ be
\begin{align*}
J_1^A&=\frac{1}{w_p}g(x_-)\\
&=\frac{1}{w_p}\left(\epsilon\log\epsilon+\theta-\epsilon(1+\log w_s)+\frac{\epsilon^2}{w_s}\right)+O(\epsilon^3)
\end{align*}

Let the $J_1$ coordinate of $B$ be $J_1^B$, and write $J_1^B=1-\epsilon^B$, for some $\epsilon^B\ll 1$. Substituting this, along with $J_2=x_-$, into the expression for the $J_1$ nullcline in~\eqref{eq:nullclines} gives
\[
x_-=\frac{1}{w_m}g(1-\epsilon^B)
\]
Expanding in terms of the small quantities $\epsilon$ and $\epsilon^B$, this gives
\begin{align*}
&\frac{\epsilon}{w_s}+O(\epsilon^2)=\\
&\frac{1}{w_m}\left(\theta-\epsilon( \log\epsilon^B +\epsilon^B+O({\epsilon^B}^2))-w_s+\epsilon^B w_s  \right)
\end{align*}
which we rearrange to find
\[
\log\epsilon^B=\frac{\theta-w_s}{\epsilon}+w_s\frac{\epsilon^B}{\epsilon}+O(1).
\]
Since we are assuming $\theta<w_s$, then $\epsilon^B$ is exponentially small, that is $\epsilon^B=O(\epsilon^n)$ for all $n\in\mathbb{N}$, thus, $J_1^B=1+O(\epsilon^n)$.

The points $A$ and $B$ collide when $J_1^A=J_1^B$, that is, when
\begin{equation} \nonumber
w_p= w_p^{SN}\equiv \theta + \epsilon \log\epsilon -\epsilon (1+\log w_s)+\frac{\epsilon^2}{w_s}+O(\epsilon^3). 
\end{equation}
\hfill $\square$

For the default parameters~\eqref{eq:defparms}, except for $w_p$, we find $w_p^{SN}=0.3027$ (4 s.f.). Note that this means for $w_p=0.3$ we are close to saddle node and there is an excitable connection with  small $\delta>0$. More generally, note that for any fixed $w_s$, as $\epsilon\rightarrow 0$ we have $w_p^{SN}\rightarrow \theta$ as expected from Theorem~\ref{thm:main}.

The following result gives an approximation of the positions of the equilibria that are created in the saddle-node bifurcation. Methods similar to those used in this proof are used in later sections for larger networks.

\begin{lemma}\label{lem:sneq}
	If $\theta<w_s$, $0<\epsilon \ll 1$, and $0<\eta\ll \frac{\epsilon}{4}$, then if  $w_p=w_p^{SN}-\eta$, the system~\eqref{eq:J2odes} has a pair of equilibria at
	\[
	(J_1,J_2)=\left(1,\frac{\epsilon}{w_s} \pm \frac{\sqrt{2\eta\epsilon}}{w_s}\right) +O(\epsilon^2).
	\]
\end{lemma}
Recall that $x_-=\frac{\epsilon}{w_s}+O(\epsilon^2)$. Thus,
\[
g''(x_-)=-\frac{w_s^2}{\epsilon}\sqrt{1-\frac{4\epsilon}{w_s}}=-\frac{w_s^2}{\epsilon}+2w_s+O(\epsilon).
\]
Using the earlier results on the location of the $J_1$ nullcline, we will have equilibria when
\[
\frac{g(J_2)}{w_p}=1+O(\epsilon^n).
\]
Expanding $g$ about $J_2=x_-$ and writing $w_p=w_p^{SN}-\eta$ gives, 
\begin{align*}
g(J_2)&=g(x_-)+\frac{(J_2-x_-)^2}{2} g''(x_-) \\
&~~~+O((J_2-x_-)^3), \\
&=w_p^{SN}-\eta+O(\epsilon^n), \\
(J_2-x_-)^2&=-\frac{2\eta}{g''(x_-)}+O(\epsilon^3),
\end{align*}
where the final line follows because $g(x_-)=w_p^{SN}$. Substituting for $g''(x_-)$ then gives the result. \hfill $\square$

\subsection{Three vertex cycle}
\label{sec:threenode}

\begin{figure}
	\begin{center}
		\includegraphics[width=3.5cm]{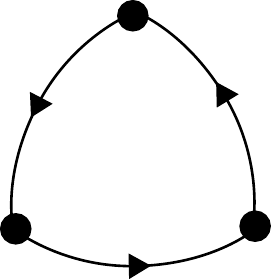}
	\end{center}
	\caption{\label{fig:GH_schematic} Graph of the three vertex cycle.}
\end{figure}

Our second example is the cycle between three vertices shown schematically in figure~\ref{fig:GH_schematic}. As a heteroclinic cycle between equilibria, this system has been studied extensively in the fields of populations dynamics~\cite{ML75}, rotating convection~\cite{BH80} and symmetric bifurcation theory~\cite{GH88}. 
We give some numerical examples of the dynamics of this system as realised by the CTRNN excitable network, and use the continuation software AUTO~\cite{auto} to show that the transition from excitable to spontaneous dynamics occurs at a saddle-node on an invariant circle (SNIC) bifurcation generating a periodic orbit.

The deterministic equations realising this graph are:
\begin{align}
\dot{y}_1 &= -y_1+ w_s \phi(y_1) + w_m \phi(y_2) +w_p \phi(y_3), \nonumber \\
\dot{y}_2 &= -y_2+ w_s \phi(y_2) + w_m \phi(y_3) +w_p \phi(y_1),  \label{eq:GHeqns} \\
\dot{y}_3 &= -y_3+ w_s \phi(y_3) + w_m \phi(y_1) +w_p \phi(y_2). \nonumber 
\end{align}
We also consider the noisy case, using the setup given in equations~\eqref{eq:ctrnn_sde}.

\begin{figure*}
	\setlength{\unitlength}{1mm}
	\begin{center}
		\begin{picture}(160,60)(0,0)
		\put(0,30){\includegraphics[trim= 2.5cm 0cm 2cm 0cm,clip=true,width=6.5cm]{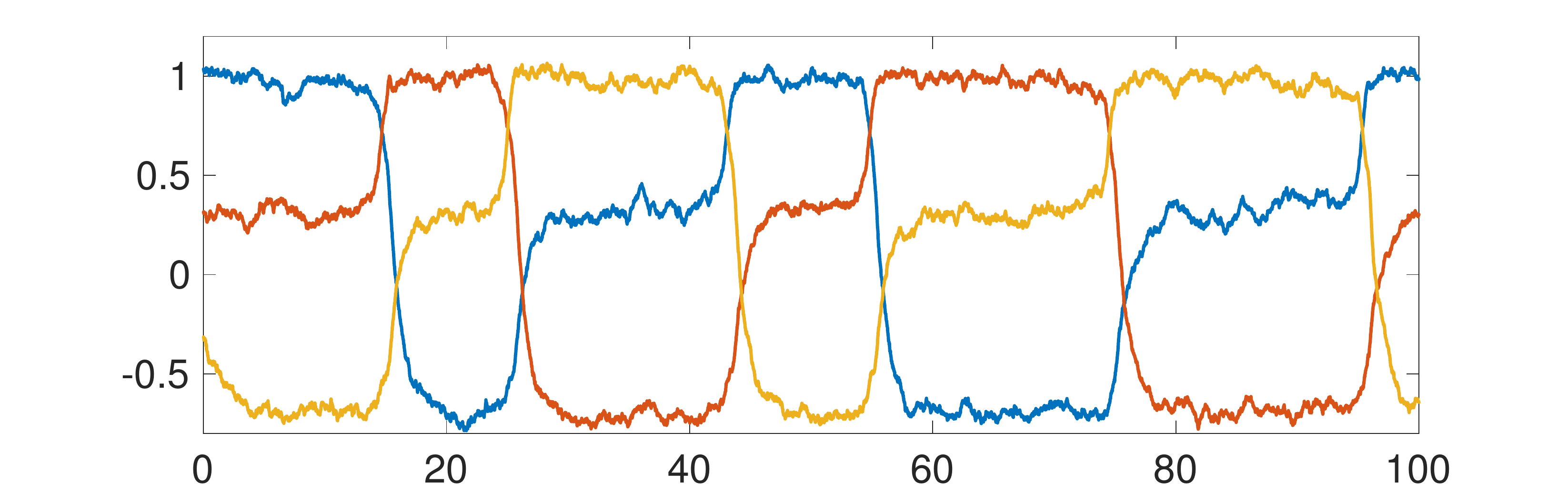}}
		\put(70,30){\includegraphics[trim= 2.5cm 0cm 2cm 0cm,clip=true,width=6.5cm]{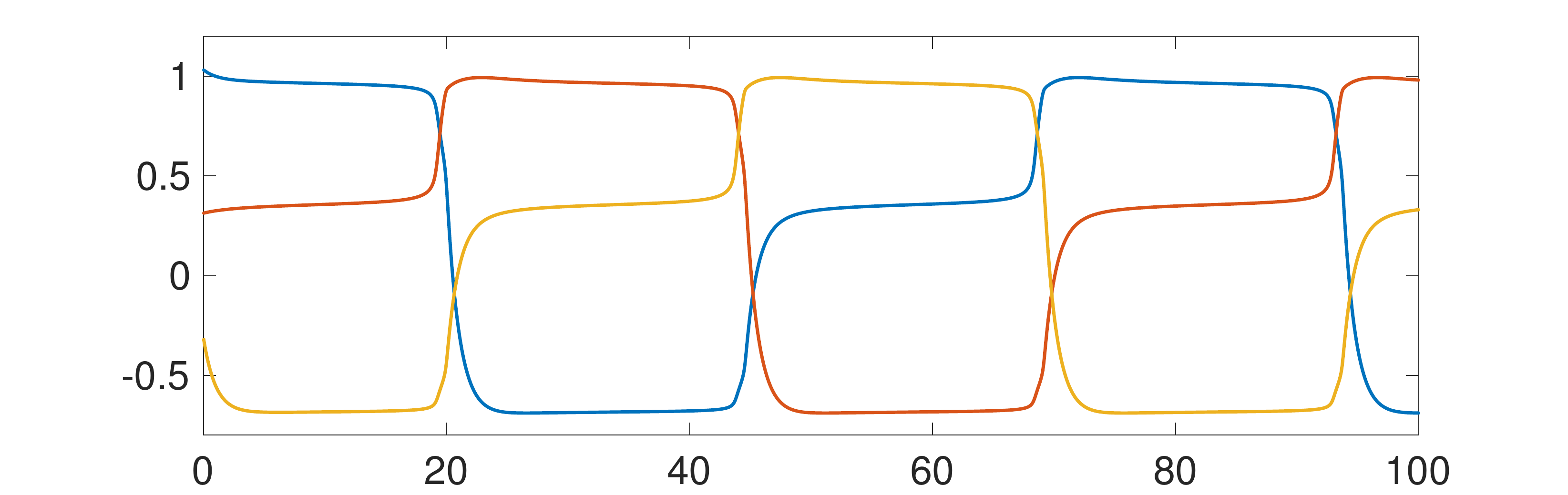}}
		\put(0,0){\includegraphics[trim= 2.5cm 0cm 2cm 0cm,clip=true,width=6.5cm]{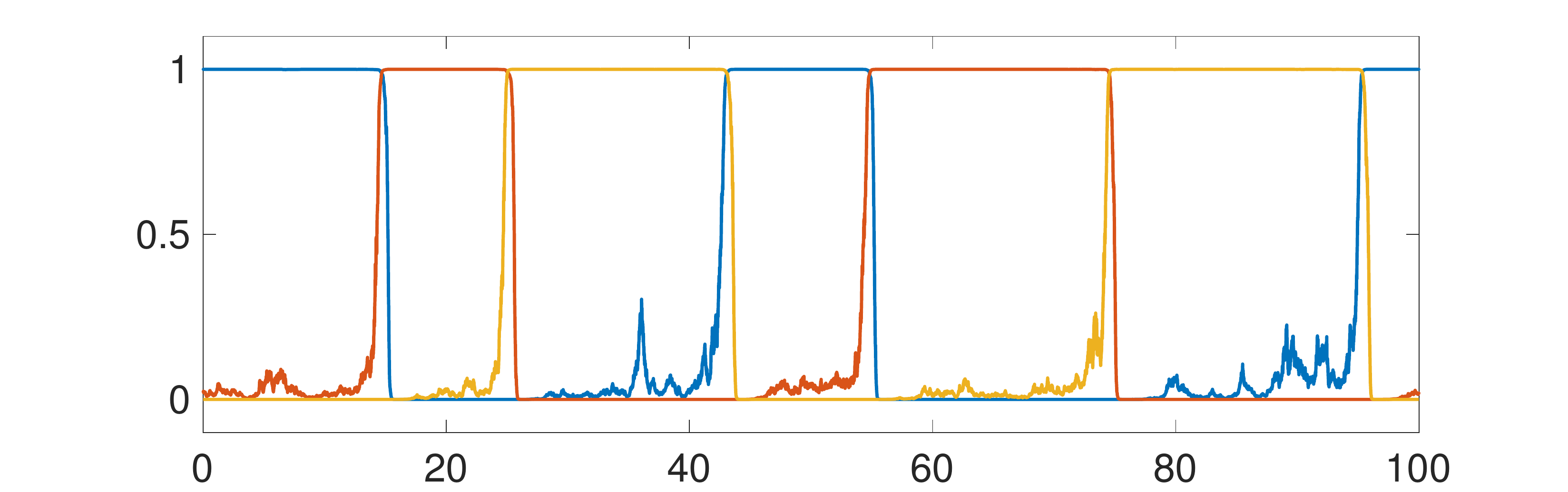}}
		\put(70,0){\includegraphics[trim= 2.5cm 0cm 2cm 0cm,clip=true,width=6.5cm]{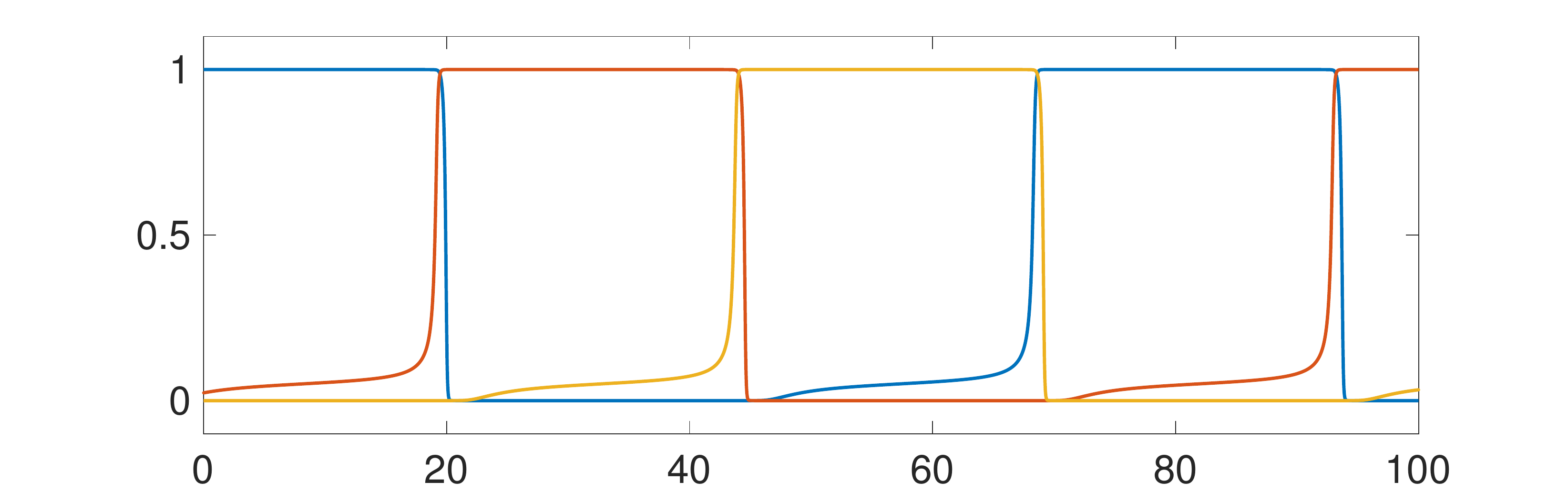}}
		
		\put(58,-1){\small $t$}
		\put(128,-1){\small $t$}
		\put(58,29){\small $t$}
		\put(128,29){\small $t$}
		
		\put(0,16){\small $J_k$}
		\put(70,16){\small $J_k$}
		\put(0,47.5){\small $y_k$}
		\put(70,47.5){\small $y_k$}

		\end{picture}
	\end{center}
	\caption{For the three-vertex cycle~\eqref{eq:GHeqns}, the top row of figures show time series of the $y_k$ variables in the excitable (left) and spontaneous (right) cases. In the left column, $w_p=0.3$ and $\sigma=0.05$. In the right column, $w_p=0.305$ and $\sigma=0$. The bottom row show the $J_k$ variables.
		Other parameters are $\epsilon=0.05$, $\theta=0.5$, $w_s=1$, $w_m=-0.7$.\label{fig:3nodets}
		}
\end{figure*}
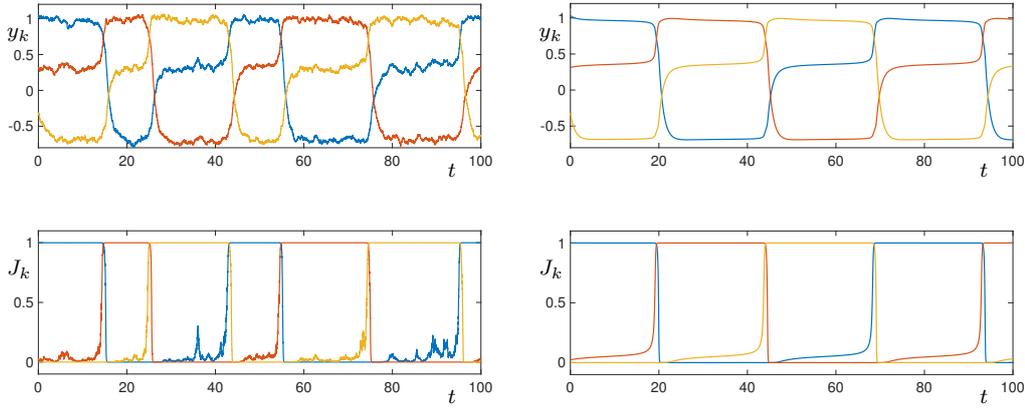

Figure~\ref{fig:3nodets} show sample time series for two different parameter sets. On the left, we show a noisy realisation with $w_p=0.3$, and on the right, a periodic solution in the deterministic system (equations~\eqref{eq:ctrnn_ode}) with $w_p=0.305$. 
Note that in both cases, for this system the $y_k$ variables oscillate between three values: high ($y_k=Y_A=w_s=1$), intermediate ($y_k=Y_L=w_p=0.3$), and low ($y_k=Y_T=w_m=-0.7$), as the cells shift between Active, Trailing and Leading. Only the first of these corresponds to $J_k\approx 1$, as can be seen in the time series plots of the $J_k$ variables in the lower panels of the figure, the other two correspond to $J_k\approx 0$.

We compute a bifurcation analysis of the system~\eqref{eq:GHeqns} using the continuation software AUTO~\cite{auto97}. 
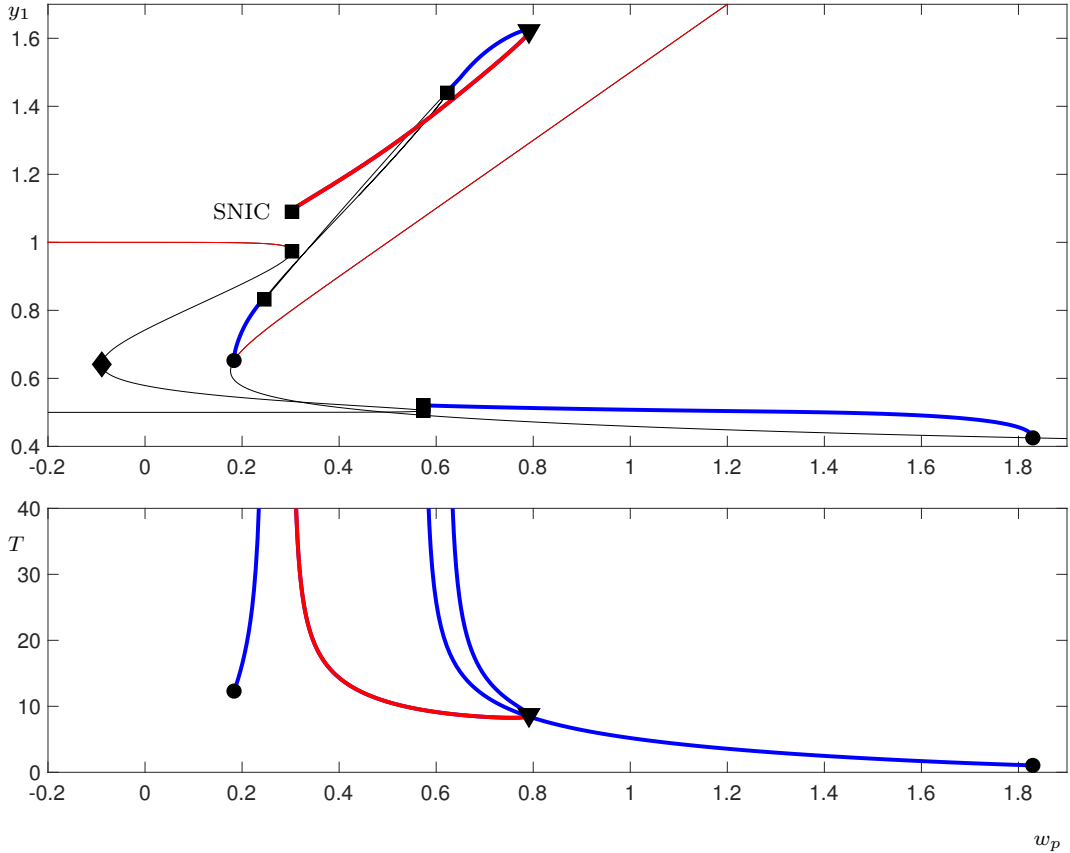
\begin{figure*}
	\setlength{\unitlength}{1mm}
	\begin{center}
		\begin{picture}(140,115)(0,0)
		\put(0,45){\includegraphics[trim= 3.5cm 0cm 3.2cm 0cm,clip=true,width=14cm]{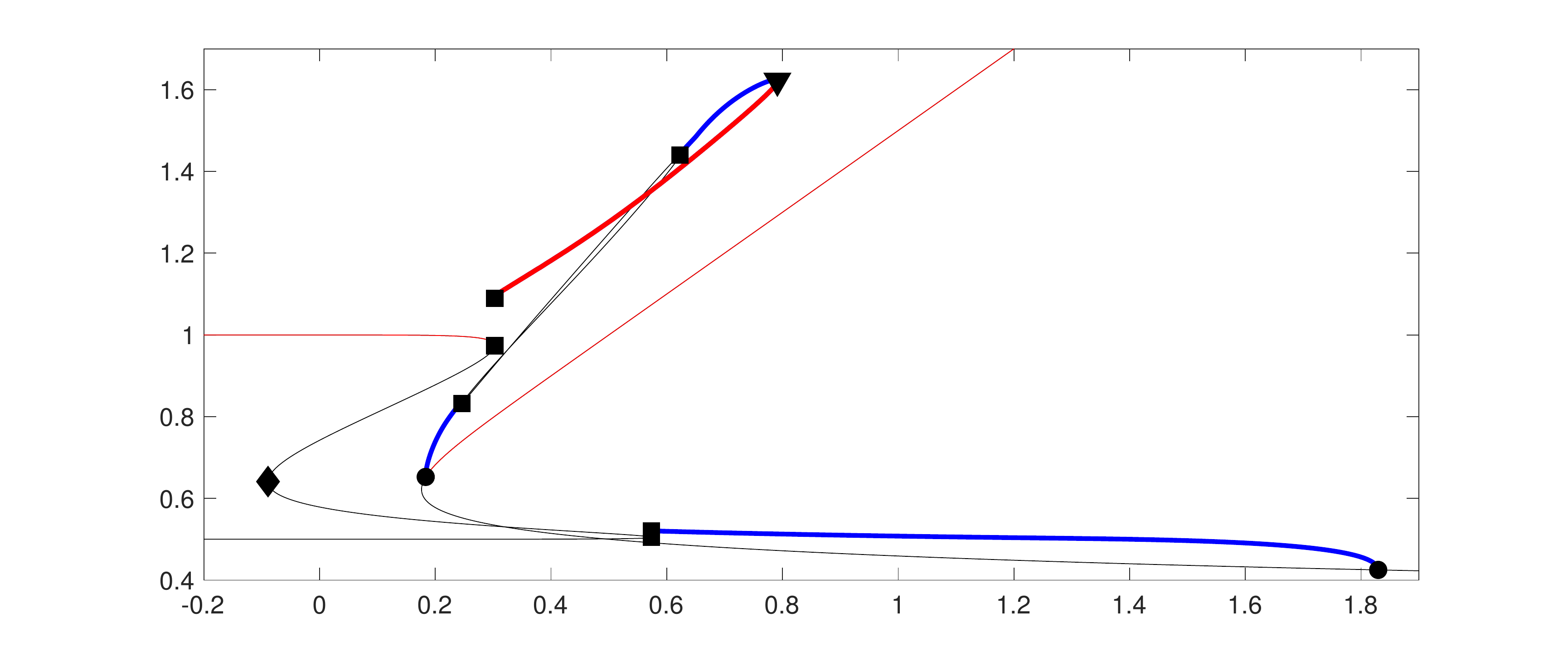}}
		\put(0,5){\includegraphics[trim= 3.5cm 0cm 3.2cm 0cm,clip=true,width=14cm]{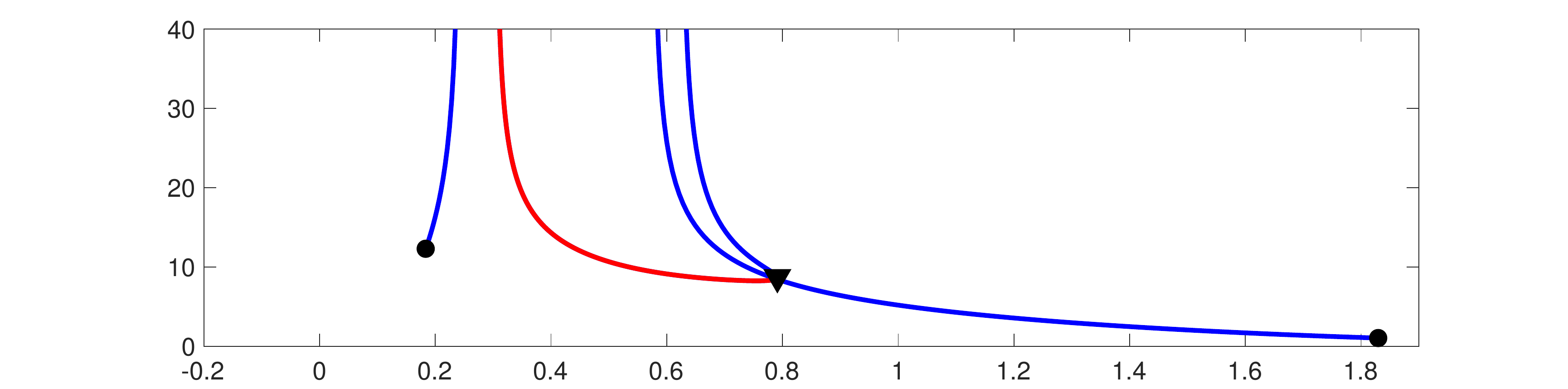}}
		
		\put(27,83){\small SNIC}
		
		\put(135,0){$w_p$}
		\put(0,39){$T$}
		\put(0,110){$y_1$}
		
		\end{picture}
		
	\end{center}
	\caption{For the three vertex cycle with equations~\eqref{eq:GHeqns}, the figure shows a bifurcation diagram as $w_p$ is varied. The top panel shows the $y_1$ coordinate of equilibrium solutions, and the maximum value of $y_1$ for periodic solutions. The lower panel shows the period of the periodic solutions. Equilibrium solutions are shown by a thin line, and periodic solutions by a thick line. Stable solutions are shown in red. Bifurcation points are indicated by various shapes: Hopf bifurcations by circles, saddle-node bifurcations by diamonds, saddle-node of periodic orbits by triangles, and saddle-node on invariant circles (SNIC) by squares. Note that there may be two squares for a single SNIC bifurcation because the maximum value of $y_1$ on the periodic orbit is not the same as the value of $y_1$ for the equilibria undergoing the saddle-node bifurcation. The SNIC bifurcation labelled at $w_p=w_p^{SNIC}\approx 0.30287$ is the transition between excitable and spontaneous dynamics.\label{fig:3nodebif}}
\end{figure*}
Figure~\ref{fig:3nodebif} shows a bifurcation diagram of this system as $w_p$ is varied. Stable solutions are shown in red. There is a saddle-node on an invariant circle (SNIC) bifurcation at $w_p=w_p^{SNIC}\approx 0.30287$. For $w_p<w_p^{SNIC}$, the diagram shows a stable equilibrium solution with $y_1\approx Y_A=1$ (and $y_2\approx Y_L$, $y_3\approx Y_T$). As $w_p$ increases through $w_p^{SNIC}$, this equilibrium disappears in a SNIC bifurcation creating a stable periodic orbit. Note that the period of the periodic orbit asymptotes to $\infty$ as the SNIC bifurcation is approached. Due to the symmetry, there are of course two further pairs of equilibria, one pair with $y_1\approx Y_T$, $y_2\approx Y_A$, $y_3\approx Y_L$, and another with $y_1\approx Y_L$, $y_2\approx Y_T$, $y_3\approx Y_A$. The symmetry causes three saddle-node bifurcations to occur simultaneously, creating the periodic orbit. If we were to instead choose the $w_p$ to be different in each of the lines in~\eqref{eq:GHeqns}, the saddle-nodes would occur independently, and a periodic orbit would exist only if all three $w_p$'s were greater than $w_p^{SNIC}$.

The time-series on the left-hand side of
figure~\ref{fig:3nodets} has $w_p=0.3<w_p^{SNIC}$. Without noise, at these parameter values, the system would remain at one equilibrium point indefinitely. The noise acts as inputs pushing the trajectory along the excitable connections. The time series on the right-hand side of figure~\ref{fig:3nodets} has $w_p=0.305>w_p^{SNIC}$, and shows the periodic orbit which has resulted from the SNIC bifurcation.

We note that this SNIC bifurcation occurs at approximately the same value of $w_p$ as the saddle-node bifurcation found in section~\ref{sec:twonode}. This is not surprising; using similar methods to those in the previous section, we can show that to lowest order in $\epsilon$, the SNIC bifurcation occurs when $w_p=w_p^{SN}$.
For $w_p<w_p^{SNIC}$ there thus exists an excitable network in the sense defined in appendix~\ref{app:networks}.

\subsection{Four node Kirk--Silber network}
\label{sec:KS}

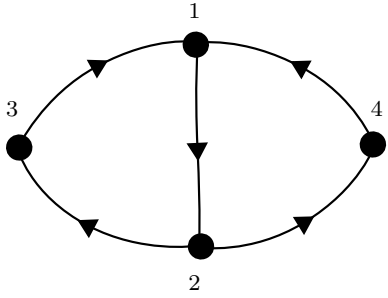
\begin{figure}
	\setlength{\unitlength}{1mm}
	\begin{center}
		\begin{picture}(50,40)(0,0)
		\put(0,4){\includegraphics[width=5cm]{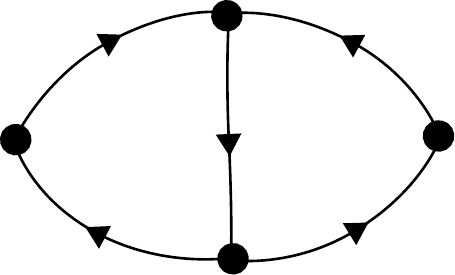}}
		\put(24,0){$2$}
		\put(24,36){$1$}
		\put(0,23){$3$}
		\put(48,23){$4$}
		\end{picture}
	\end{center}
	\caption{\label{fig:KS_schematic} Graph of the four vertex Kirk--Silber network.}
\end{figure}

For our next example, we consider a graph with the structure of the  Kirk--Silber network~\cite{KS94}, shown schematically in figure~\ref{fig:KS_schematic}. This graph has one vertex which has two outgoing edges, and the dynamics here are somewhat different to vertices with only one outgoing edge. The bulk of this section is devoted to an analysis of these differences, in particular, the possibility of an additional equilibrium in the network attractor with two active cells.

The corresponding deterministic equations for this network are (moving immediately into the $J_i$ variables):
\begin{align}
\dot{J}_1 &=\frac{1}{\epsilon}J_1(1-J_1)\left(-\phi^{-1} J_1+  w_s J_1  + w_m J_2  +w_p J_3 +w_p J_4\right), \label{eq:KSJ} \\
\dot{J}_2 &=\frac{1}{\epsilon}J_2(1-J_2)\left( -\phi^{-1} J_2 + w_p J_1 + w_s J_2 +w_m J_3+w_m J_4 \right), \nonumber \\
\dot{J}_3 &= \frac{1}{\epsilon}J_3(1-J_3)\left( -\phi^{-1} J_3 + w_m J_1 + w_{p_3} J_2 +w_s J_3 +w_t J_4 \right), \nonumber \\
\dot{J}_4 &= \frac{1}{\epsilon}J_4(1-J_4)\left( -\phi^{-1} J_4 + w_m J_1 + w_{p_4} J_2 +w_t J_3 +w_s J_4\right). \nonumber
\end{align}

We can break the symmetry between $J_3$ and $J_4$ by choosing $w_{p_3}\neq w_{p_4}$. In fact, in what follows, we will frequently write  $w_{p_3}=w_{p_4}+\Delta w$, for $\Delta w>0$, and choose $w_{p_4}=w_p$ for simplicity.

We consider first the dynamics near each of the vertices that have exactly one outgoing edge (vertices 1, 3 and 4 in the graph; see figure~\ref{fig:KS_schematic}). Again using the same techniques that were used in section~\ref{sec:twonode} (lemma~\ref{lem:sneq}), we can show that for  $w_p,w_m,w_t<\theta<w_s$, and $w_p=w_p^{SN}+\eta$, $0<\eta<\frac{\epsilon}{4}$, there exist equilibria solutions at, for example
\[
(J_1,J_2,J_3,J_4)=\left(1,\frac{\epsilon}{w_s}\pm\frac{\sqrt{2\eta\epsilon}}{w_s},0,0\right) +O(\epsilon^2).
\]
That is, there is a transition from excitable to spontaneous dynamics (in this case between cells $1$ and $2$, but also between cells $3$ and $1$ and cells $4$ and $1$) as $w_p$ is increased through $w_p^{SN}$.

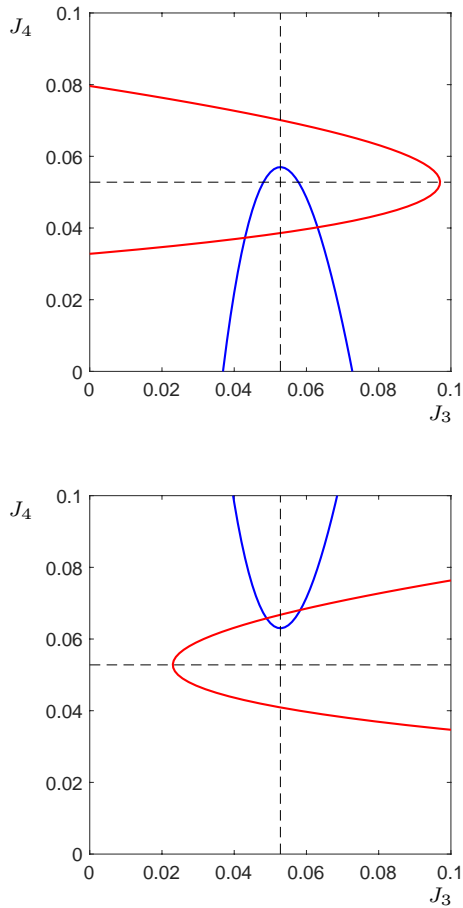
\begin{figure}
	\setlength{\unitlength}{0.8mm}
	\begin{center}
		\begin{picture}(80,155)(0,5)
		\put(5,5){\includegraphics[trim= 2.5cm 0cm 2cm 0cm,clip=true,width=6cm]{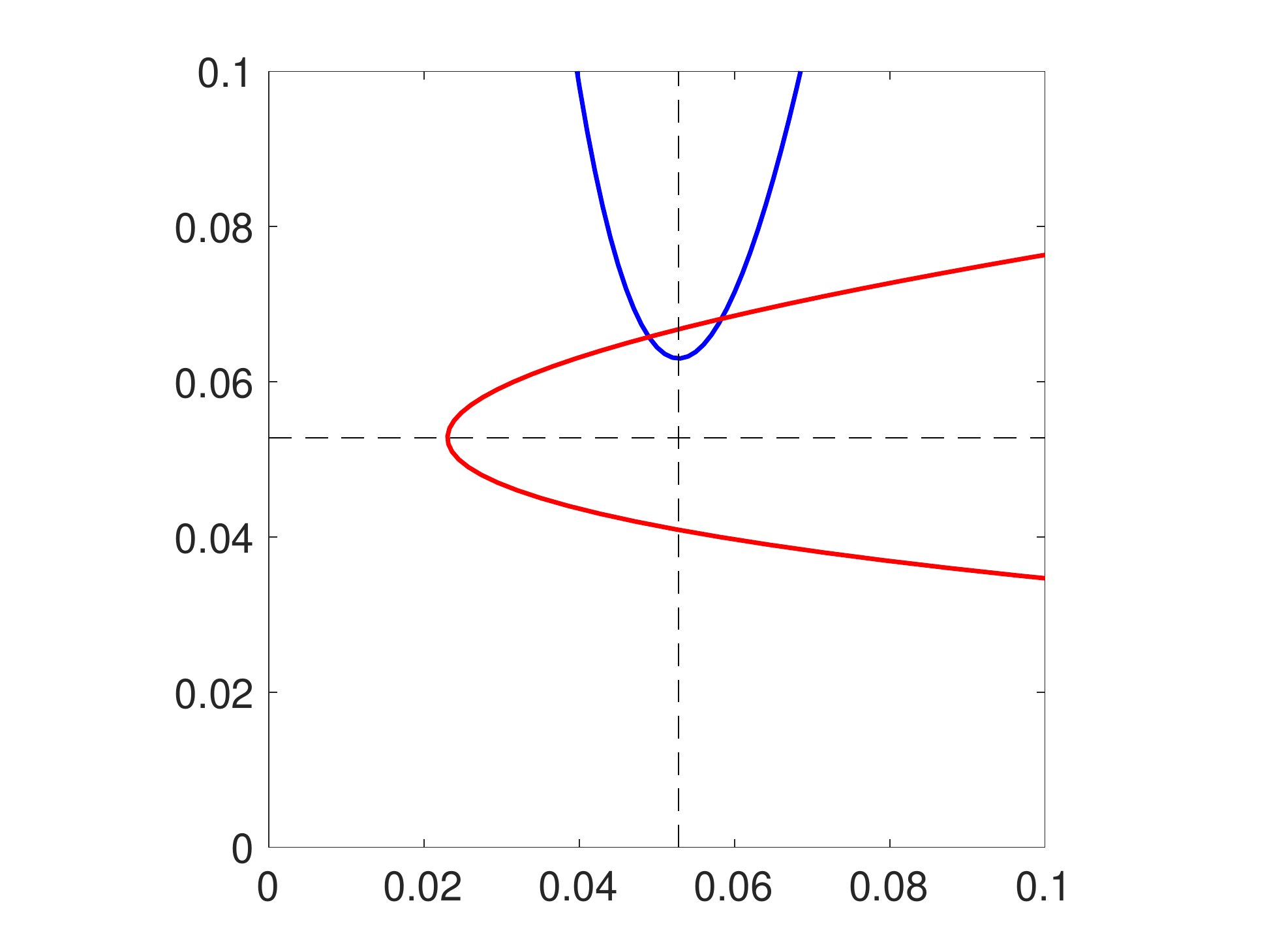}}
		\put(5,85){\includegraphics[trim= 2.5cm 0cm 2cm 0cm,clip=true,width=6cm]{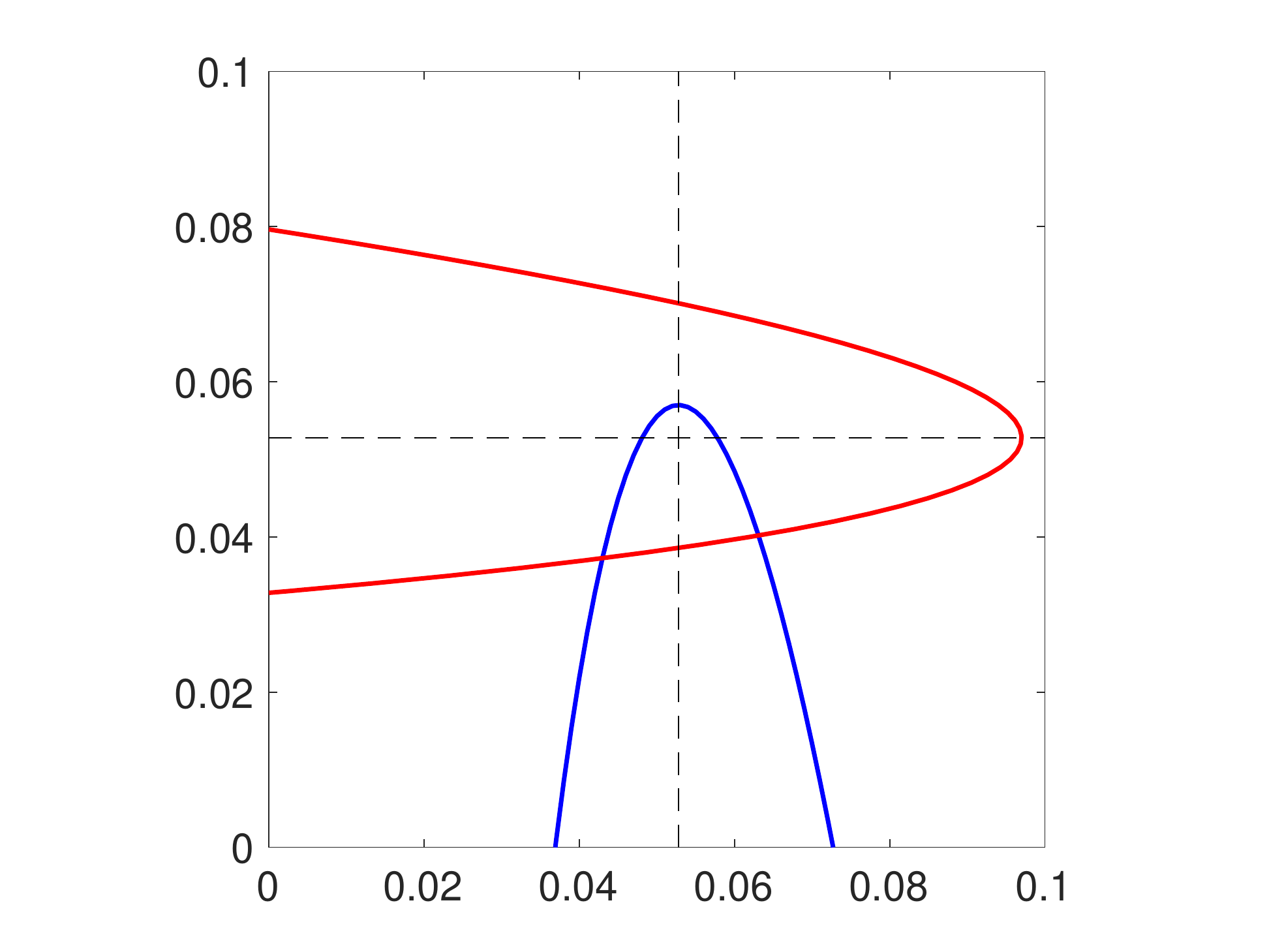}}
		\put(69,5){$J_3$}
		\put(69,85){$J_3$}
		\put(0,69){$J_4$}
		\put(0,149){$J_4$}
		%\put(0,12){$x_-$}
		%\put(64,15){$A$}
		%\put(69.5,15){$B$}
		\end{picture}
		
	\end{center}
	\caption{\label{fig:j34null} Schematic illustration of the nullclines of the system~\eqref{eq:KSJ234} system along the surface $J_2=1$.  The blue curves are the $J_3$ nullclines, and the red curves are the $J_4$ nullclines.  In the upper panel, $w_t=0.05$, $w_{p_3}=0.30$, $w_{p_4}=0.298$. In the lower panel, $w_t=-0.05$, $w_{p_3}=0.306$, $w_{p_4}=0.304$. The dashed lines show the where the nullcline would lie if $w_t=0$.  }
\end{figure}

The dynamics close to the vertex with two outgoing connections (vertex 2) is modified by the presence of the additional parameter $w_t$. Consider the three dimensional subsystem of~\eqref{eq:KSJ} with $J_1=0$, that is:
\begin{align}
\dot{J}_2 &=\frac{1}{\epsilon}J_2(1-J_2)\left( -\phi^{-1} J_2 + w_s J_2 +w_m J_3+w_m J_4 \right) \label{eq:KSJ234} \\
\dot{J}_3 &= \frac{1}{\epsilon}J_3(1-J_3)\left( -\phi^{-1} J_3 + w_{p_3} J_2 +w_s J_3 +w_t J_4 \right) \nonumber \\
\dot{J}_4 &= \frac{1}{\epsilon}J_4(1-J_4)\left( -\phi^{-1} J_4 + w_{p_4} J_2 +w_t J_3 +w_s J_4\right). \nonumber
\end{align}
We can perform a similar calculation to that shown in section~\ref{sec:twonode} to show that 
there is a section of the $J_2$ null-surface which lies asymptotically close to the surface $J_2=1$. Equilibria solutions exist on this null-surface if the $J_3$ and $J_4$ null-surfaces intersect there, that is, if there are solutions to the pair of equations
\begin{align}
g(J_3)&=w_{p_3}+w_tJ_4, \\
g(J_4)&=w_{p_4}+w_tJ_3. 
\end{align}
We assume without loss of generality that $w_{p_3}>w_{p_4}$ (i.e.~$\Delta w>0$) and then the arrangement of these curves is in one of the configurations shown in figure~\ref{fig:j34null}. If $w_t<0$ (lower panel), equilibria solutions exist (i.e.~the red and blue curves intersect) for a range of $w_{p_3}$ and $w_{p_4}$ with both larger than $w_p^{SN}$: that is, the transition to spontaneous dynamics happens at a larger value of $w_{p_j}$ (than if $w_t=0$). If $w_t>0$ (upper panel), the opposite happens: that is, the transition to spontaneous dynamics occurs at a smaller value of $w_{p_j}$. Solving these equations exactly requires solving a quartic equation, and the resulting expression is not illuminating. We label the value of $w_p$ at which this transition from spontaneous to excitable dynamics occurs as $w_p^{SN'}$, and note that this is a function of $w_t$, $w_s$, $\epsilon$, $\theta$ as well as more generally, the number of Leading directions from that cell.

For the specific system~\eqref{eq:KSJ}, with $w_{p_3}=w_{p_4}+\Delta w=w_p+\Delta w$, we thus have two conditions. If $w_p<\min(w_p^{SN},w_p^{SN'}-\Delta w)$ then the system is excitable along all connections. If $w_p>\max(w_p^{SN},w_p^{SN'}-\Delta w)$ then a periodic orbit will exist. If neither of these conditions holds then we will see excitable connections in some places and spontaneous transitions in others. 

\begin{figure*}
	\setlength{\unitlength}{1mm}
	\begin{center}
		\begin{picture}(160,60)(0,0)
		\put(0,30){\includegraphics[trim= 2.5cm 0cm 2cm 0cm,clip=true,width=6.5cm]{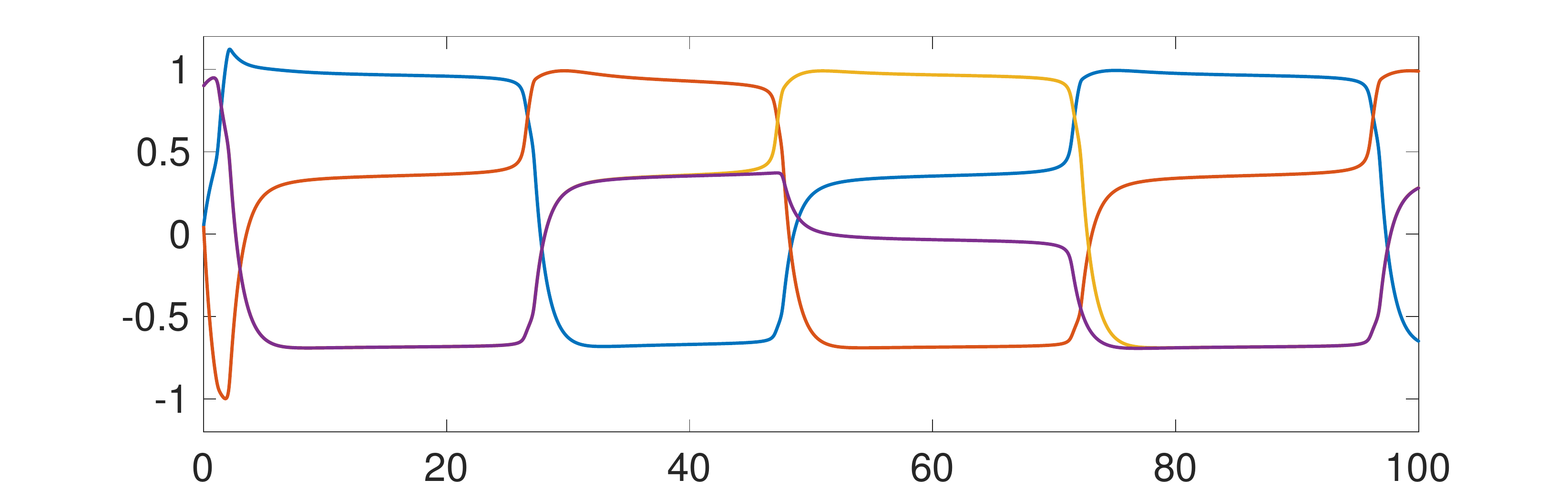}}
		\put(70,30){\includegraphics[trim= 2.5cm 0cm 2cm 0cm,clip=true,width=6.5cm]{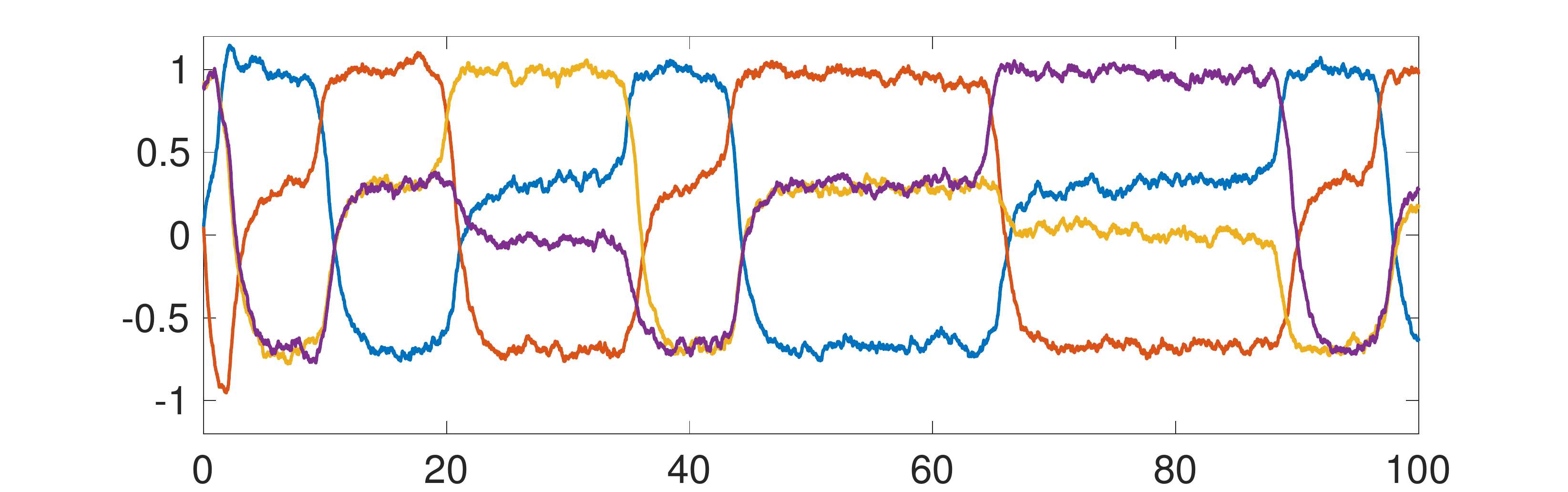}}
		\put(0,0){\includegraphics[trim= 2.5cm 0cm 2cm 0cm,clip=true,width=6.5cm]{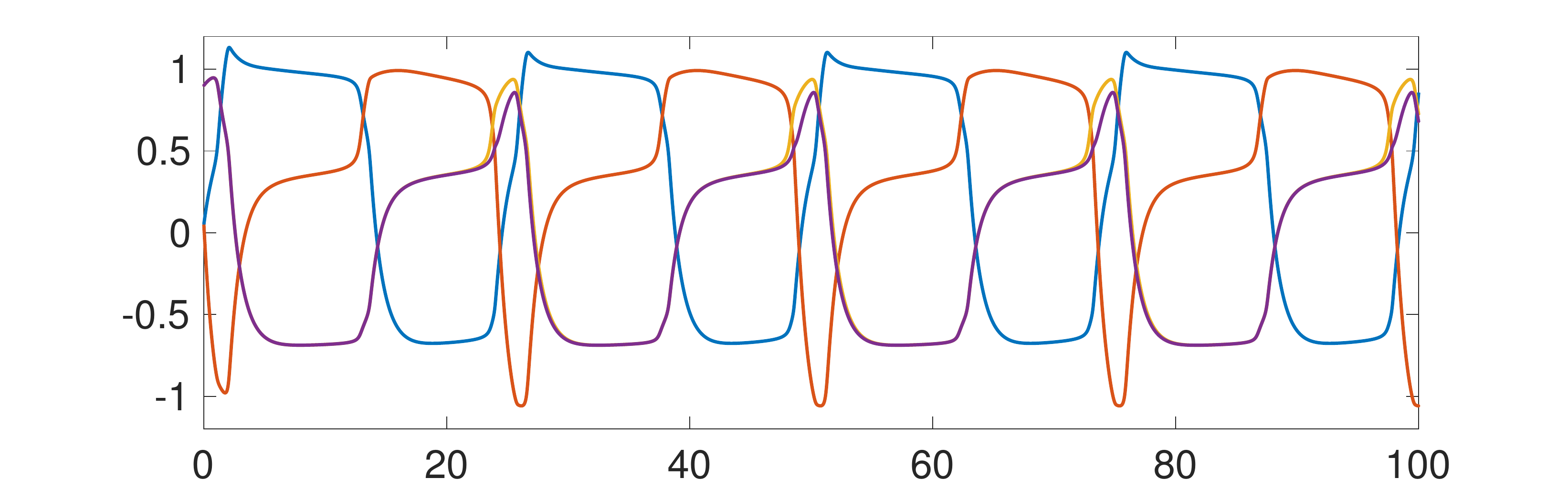}}
		\put(70,0){\includegraphics[trim= 2.5cm 0cm 2cm 0cm,clip=true,width=6.5cm]{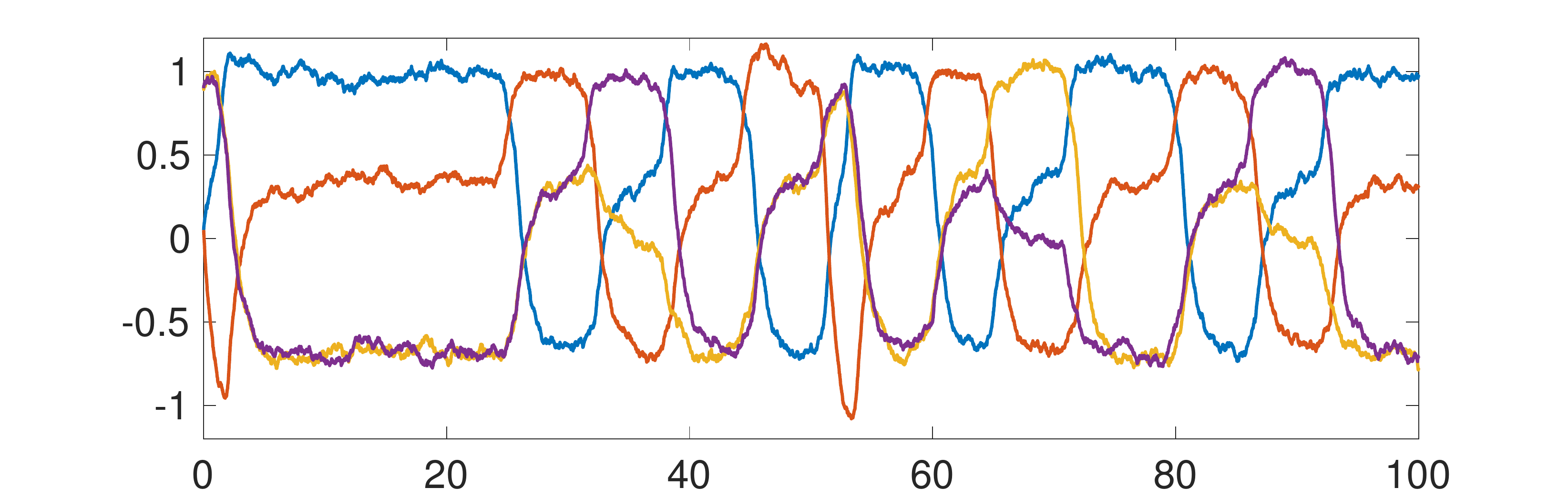}}
		
		\put(0,55){(a)}
		\put(70,55){(b)}
		\put(0,25){(c)}
		\put(70,25){(d)}
		
		\put(58,-1){\small $t$}
		\put(128,-1){\small $t$}
		\put(58,29){\small $t$}
		\put(128,29){\small $t$}
		
		\put(-2,18){\small $y_j$}
		\put(68,18){\small $y_j$}
		\put(-2,48){\small $y_j$}
		\put(68,48){\small $y_j$}
		
		\end{picture}
	\end{center}
	\caption{The figures show time series of simulations of the $y_j$ variables in the Kirk--Silber type network~\eqref{eq:KSJ}. In (a), $w_p=0.305$, and $\sigma=0$; in (b), $w_p=0.3$ and $\sigma=0.05$; in (c) $w_p=0.315$, and $\sigma=0$; in (d), $w_p=0.315$ and $\sigma=0.05$.\label{fig:KSts}
	}
\end{figure*}

In figure~\ref{fig:KSts} we show some example time-series of the system~\eqref{eq:KSJ} (in the $y_k$ coordinates). In panel (a), parameters are such that $w_p>\max(w_p^{SN},w_p^{SN'}-\Delta w)$, so we see a periodic solution in the deterministic system. Note that the $y_3$ (yellow) coordinate becomes close to $Y_A=1$ during this trajectory, but the $y_4$ (purple) coordinate does not: it switches between $Y_L=0.3$, $Y_D=0$ and $Y_T=0.7$. In panel (b), parameters are such that  $w_p<\min(w_p^{SN},w_p^{SN'}-\Delta w )$, so without noise the trajectory would remain at a single equilibrium solution. Here we add noise with $\sigma=0.05$, and the trajectory can be seen exploring the network. Note that there are some transitions between $\xi_2$ ($y_2$ is red) and $\xi_3$ ($y_3$ is yellow), and some from $\xi_2$ to $\xi_4$ ($y_4$ is purple). In panels (c) and (d), we increase $w_p$ further away from the saddle-node bifurcation (further into the regime of spontaneous transitions) and observe some qualitative differences in the trajectories. In the deterministic case (c), the periodic solution now transitions from near the bottleneck region $P_2$ to a region of phase space where $y_3$ and $y_4$ are both Active. We label this region of phase space as $P_{3,4}$. In the noisy case (d) (which is also in the spontaneous regime), the trajectory makes transitions from the bottleneck region $P_2$ to each of $P_3$, $P_4$, and to  $P_{3,4}$.

\begin{figure}
	\setlength{\unitlength}{1mm}
	\begin{center}
		\begin{picture}(80,50)(0,0)
		\put(0,28){\includegraphics[trim= 2.5cm 0cm 2cm 0cm,clip=true,width=8cm]{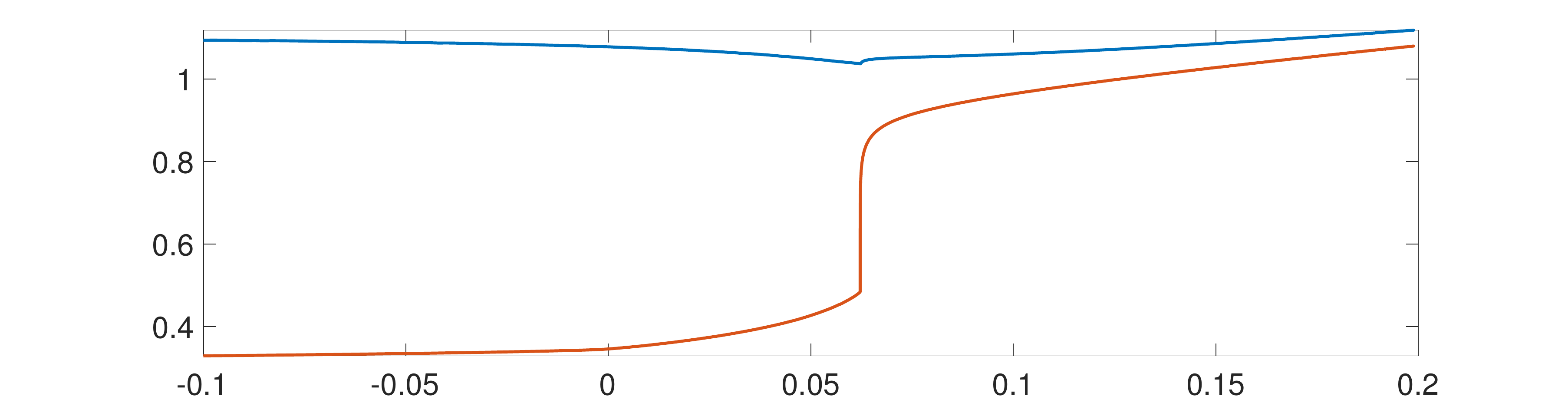}}
		\put(0,0){\includegraphics[trim= 2.5cm 0cm 2cm 0cm,clip=true,width=8cm]{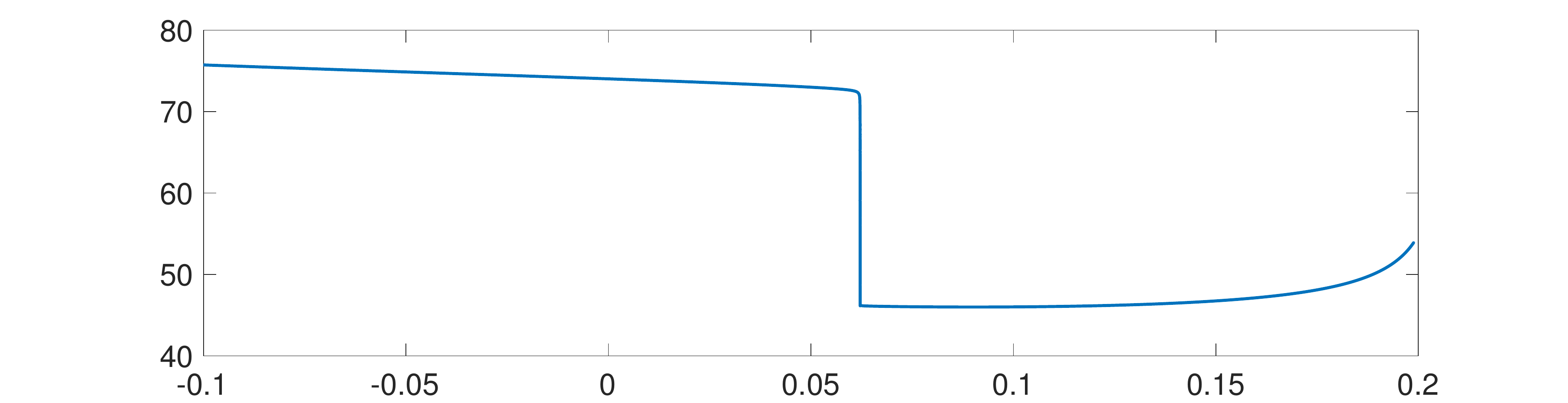}}
		\put(0,20){$T$}
		
		\put(0,31){\rotatebox{90}{$\max y_3$, $y_4$}}
	
		\put(70,0){$w_t$}
		\put(70,28){$w_t$}
		\end{picture}
	\end{center}
	\caption{\label{fig:KSt_bif}Behaviour of periodic orbits in the  Kirk--Silber-type four-node network~\eqref{eq:KSJ} as the parameter $w_t$ is varied. The top panel shows $\max(y_3)$ (blue) and $\max(y_4)$ (red) along the periodic orbit. The bottom panel shows the period of the orbit on varying $w_t$. 
	}
\end{figure}
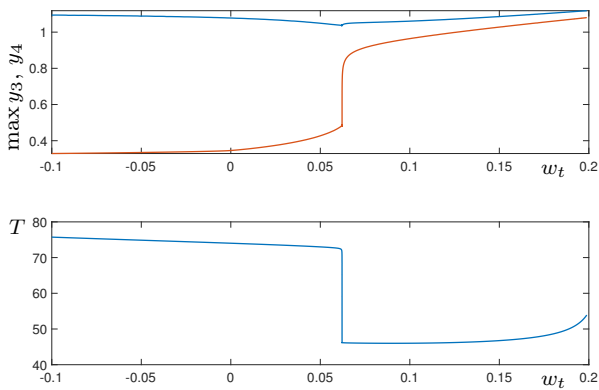

 In the above simulations, we have used $w_t=0$, but we observe that the type of qualitative behaviour observed depends both on the parameters $w_p$ and $w_t$. Specifically, a sufficiently negative $w_t$ provides a suppression effect, meaning that only a single cell $y_k$ can be active at any one time, but the transitional value of $w_t$ depends on $w_p$.
In figure~\ref{fig:KSt_bif} we show maximum values of $y_3$ and $y_4$ along the periodic orbits as $w_t$ is varied. It can be seen clearly here that the transition between periodic orbits which visit $P_3$ (where $\max(y_3)$ is significantly larger than $\max(y_4)$) and those which visit $P_{3,4}$ (where $\max(y_3)\approx\max(y_4)$) is quite sharp. 

\begin{figure*}
	\setlength{\unitlength}{1mm}
	\begin{center}
		\begin{picture}(150,60)(0,0)
		\put(0,0){\includegraphics[trim= 2.cm 0cm 2.5cm 0cm,clip=true,width=15cm]{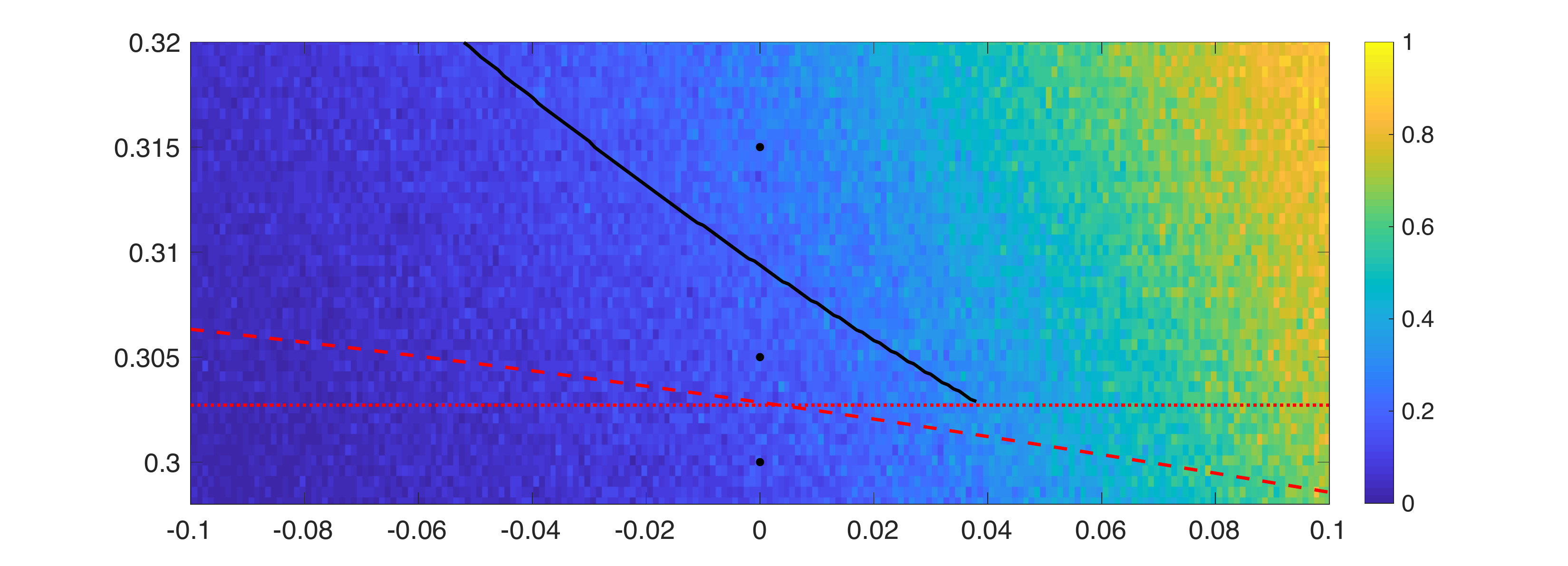}}
		\put(130,1){$w_t$}
		\put(0,53){$w_p$}
		\put(67,11){\color{white}{(b)}}
		\put(67,23){\color{white}{(a)}}
		\put(75,46){(c),(d)}
		\end{picture}
	\end{center}
	\caption{ \label{fig:KSt} Behaviour of the Kirk--Silber-type four-node network~\eqref{eq:KSJ} as the parameters $w_t$ and $w_p$ are varied. The red lines are the curves $w_p=w_p^{SN}$ (dotted) and $w_p=w_p^{SN'}$ (dashed; determined numerically by solving a quartic equation). For $w_p$ above both of these lines periodic solutions exist in the deterministic system. To the left of the black line these periodic solutions visit $P_3$, to the right they visit $P_{3,4}$. The background colours are results from noisy simulations with $\sigma=0.05$. The colour indicates the ratio of transitions to $P_{3,4}$ to the total number of transitions to $P_3$, $P_4$ and $P_{3,4}$. The labelled dots give the parameter values of the time-series plots in figure~\ref{fig:KSts}.}
\end{figure*}

We extend these results to show the behaviour as both parameters $w_p$ and $w_t$ are varied in figure~\ref{fig:KSt}. The data in this figure shows the observed behaviours for both noisy and deterministic systems as the parameters $w_p$ and $w_t$ are varied. The red lines are the curves $w_p=w_p^{SN}$ (dotted) and $w_p=w_p^{SN'}$ (dashed). If $w_p$ is above both of these lines then all transitions are spontaneous, and so a  periodic orbit exists in the system. If $w_p$ lies below either one (or both) of these lines, then at least one of the transitions will be excitable and so there will be no periodic solutions. The black line shows the boundary between those periodic solutions which visit $P_3$ (to the left of the black line) and those which visit $P_{3,4}$ (to the right of the black line), as determined by the location of the sharp transition in calculations similar to those shown in figure~\ref{fig:KSt_bif} for a range of $w_p$. The background colours are results from noisy simulations. The colour indicates the ratio of transitions to $P_{3,4}$ to the total number of transitions to $P_3$, $P_4$ and $P_{3,4}$. Interestingly, the noisy solutions require a much larger value of $w_t$ than the deterministic ones to have a significant proportion of transitions to $P_{3,4}$.

\begin{figure*}
	\setlength{\unitlength}{0.8mm}
	\begin{center}
		\begin{picture}(160,70)(0,0)
		\put(4,0){\includegraphics[trim= 2.5cm 1cm 1cm 1cm,clip=true,width=5.7cm]{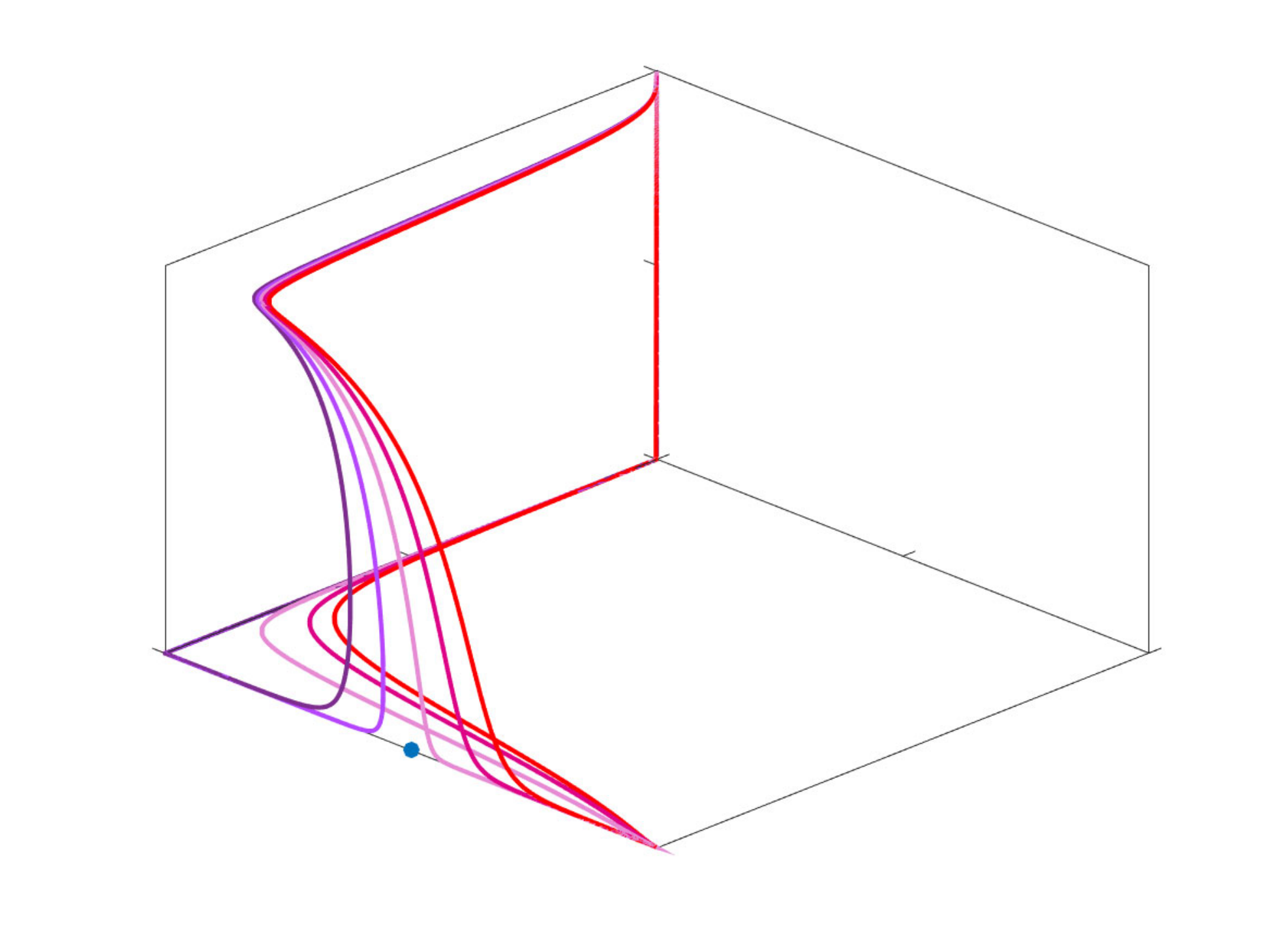}}
		\put(86,0){\includegraphics[trim= 2.5cm 1cm 1cm 1cm,clip=true,width=5.7cm]{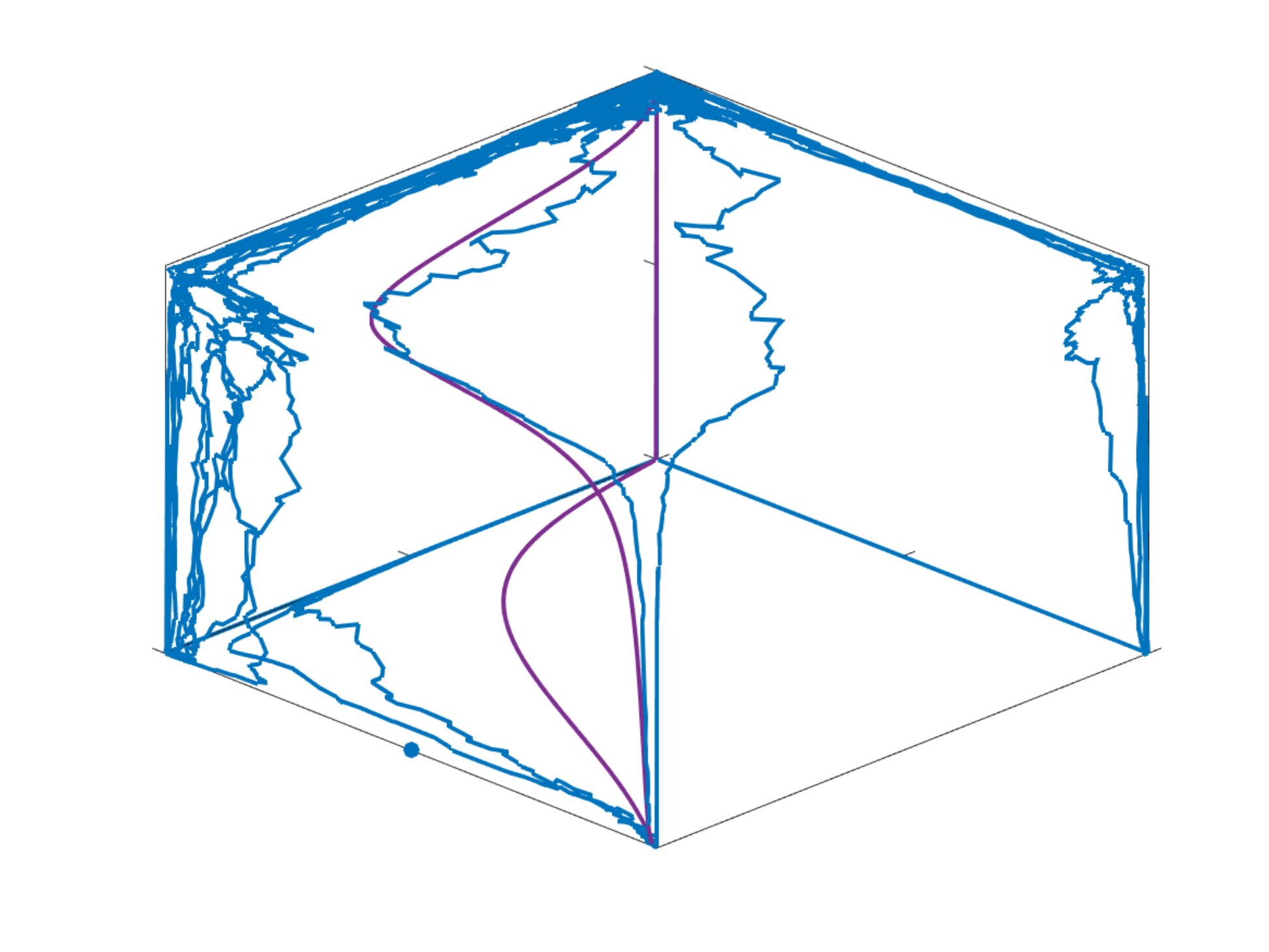}}
		\put(0,13){$J_3$}
		\put(34,59){$J_2$}
		\put(74,14){$J_4$}
		
		\put(82,13){$J_3$}
		\put(116,59){$J_2$}
		\put(156,14){$J_4$}
		
		\put(0,59){(a)}
		\put(82,59){(b)}
		
		\put(30,45){\vector(-3,-1){10}}
		\put(40,30){\vector(0,1){10}}
		\put(25,22){\vector(3,1){10}}
		
		\put(110,54){\vector(-3,-1){10}}
		\put(130,54){\vector(2,-1){10}}
		
		\put(85,33){\vector(0,-1){10}}
		\put(155,33){\vector(0,-1){10}}
		
		\put(112,10){\vector(-2,1){10}}
		\put(140,18){\vector(-2,1){10}}
		
		%\put(160,70){x}
		\end{picture}
	\end{center}
	\caption{ \label{fig:KSJ234} The figures show trajectories for the system~\eqref{eq:KSJ} projected onto the three-dimensional space with $J_1=0$. In panel (a), five periodic trajectories in the deterministic system are shown for different values of $w_p$, ranging from $w_p=0.309$ (left curve, dark purple), to $w_p=0.3096$ (right curve, red) in increments of $0.002$.  In panel (b), we set $w_p=0.315$ (as in panels (c) and (d) of figure~\ref{fig:KSts}). Noisy trajectories ($\sigma=0.05$) are shown in blue, and the periodic orbit of the deterministic system is shown in purple. In both panels, an equilibrium of the three-dimensional system~\eqref{eq:KSJ234} is shown by a blue dot. Other parameters are $\eps=0.05$, $\theta=0.5$, $w_t=0$, $w_s=1$, $w_m=-0.5$. The arrows indicate the direction of flow.}
	\end{figure*}
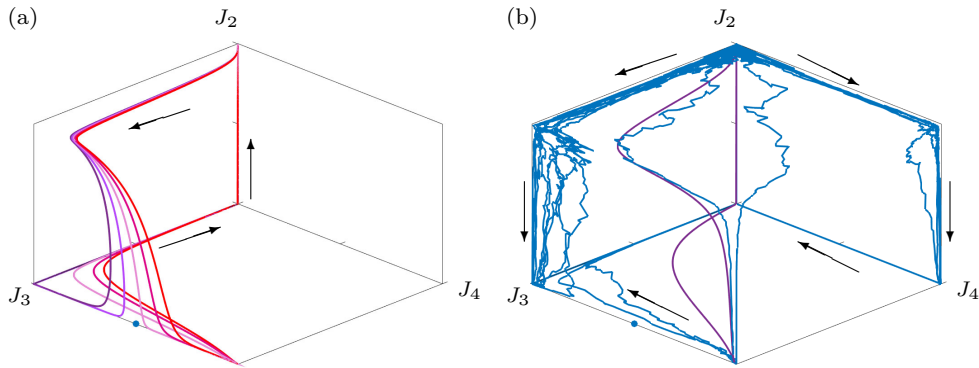

These changes in qualitative dynamics can be explained in terms of the three-dimensional subsystem with $J_1=0$, given by equations~\eqref{eq:KSJ234}.
In this three-dimensional system, there are stable equilibria at 
\begin{align*}
(J_2,J_3,J_4)&=(0,1,0)+O(\epsilon)\\ (J_2,J_3,J_4)&=(0,0,1)+O(\epsilon)\\
(J_2,J_3,J_4)&=(0,1,1)+O(\epsilon)
\end{align*}
as well as further unstable/saddle equilibria. Recall that these equilibria are not on the boundaries of the box (which are not part of the domain). In figure~\ref{fig:KSJ234} we show solutions from the full four-dimensional system~\eqref{eq:KSJ} projected onto the three-dimensional space with $J_1=0$. In panel (a), we show the periodic solutions from the deterministic systems for five different values of $w_p$, ranging from $w_p=0.309$ (left curve, dark purple), to $w_p=0.3096$ (right curve, red) in increments of $0.0002$. It can be seen that the first two of these trajectories approach the saddle equilibria (marked as a blue dot) from one side of its stable manifold, and the latter three from the other. The first three thus visit $P_3$, and the latter three visit $P_{3,4}$. It is the transition of the periodic orbit across the stable manifold which results in the rapid change in the qualitative behaviour of the periodic orbit, and likely indicates that a homoclinic bifurcation to this saddle point separates these behaviours. That is, the sharp transition in $T$ in figure~\ref{fig:KSt_bif} should actually extend to $\infty$ on both sides.  In panel (b), we show both noisy and deterministic trajectories with $w_p=0.315$. Note that only one of the noisy trajectories follows the deterministic trajectory closely: the majority visit $P_3$.

\subsection{A ten node network}
\label{sec:tennode}

\begin{figure}
	\begin{center}
		{\includegraphics[width=7cm]{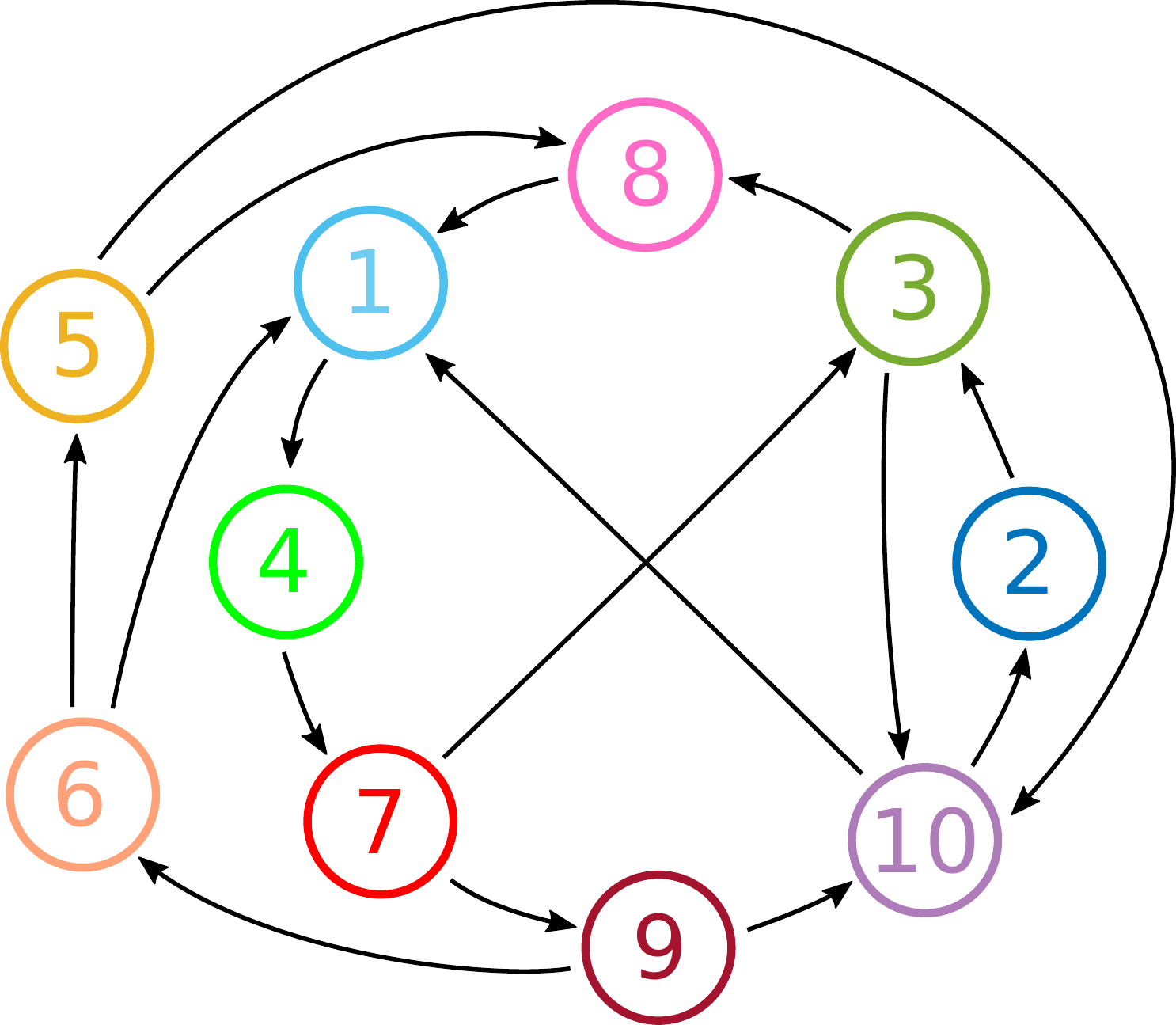}}
	\end{center}
	\caption{A example of a directed graph between ten nodes with no loops of order one or two and no $\Delta$-cliques. Three timeseries from  realisations of this as a network using CTRNNs are shown in  figures~\ref{fig:tennodetraj_det}-\ref{fig:tennodetraj_noisy2}.
	\label{fig:tennodegraph}}
\end{figure}

In this section, we demonstrate the method of construction described in section~\ref{sec:method} for a larger network. Specifically, we randomly generated a directed graph between 10 vertices, with the constraints that it contained no one-loops, two-loops or $\Delta$-cliques, and such that the graph does not have feed-forward structure (i.e. you cannot get `stuck' in a subgraph by following the arrows). The graph we consider is shown in figure~\ref{fig:tennodegraph}.

We ran one simulation of the deterministic CTRNN system~\eqref{eq:ctrnn_ode}, and two simulations of the noisy CTRNN system~\eqref{eq:ctrnn_sde}, and in each case, randomly generated the entries for the $w_p$ in the equation~\eqref{eq:w}. For the deterministic system, the entries of $w_p^{ij}$ were chosen independently from the uniform distribution $U(0.32,0.34)$. For the noisy systems, the entries of $w_p^{ij}$ were chosen independently from the uniform distribution $U(0.30,0.32)$. The remaining parameters were set at the default parameter values given in~\eqref{eq:defparms} except $w_t=-0.3$ for the deterministic system, and one of the noisy systems, and $\sigma=0.01$ for the noisy systems. Note that for the deterministic parameter values there are bottlenecks in the phase space close to the locations in phase space for the excitable states that are present for the parameter values used in the stochastic cases. 

\begin{figure*}
	\setlength{\unitlength}{1mm}
	\begin{center}
		\begin{picture}(140,85)(0,5)
		\put(0,5){\includegraphics[trim= 3.5cm 0.5cm 3.2cm 0.5cm,clip=true,width=14cm]{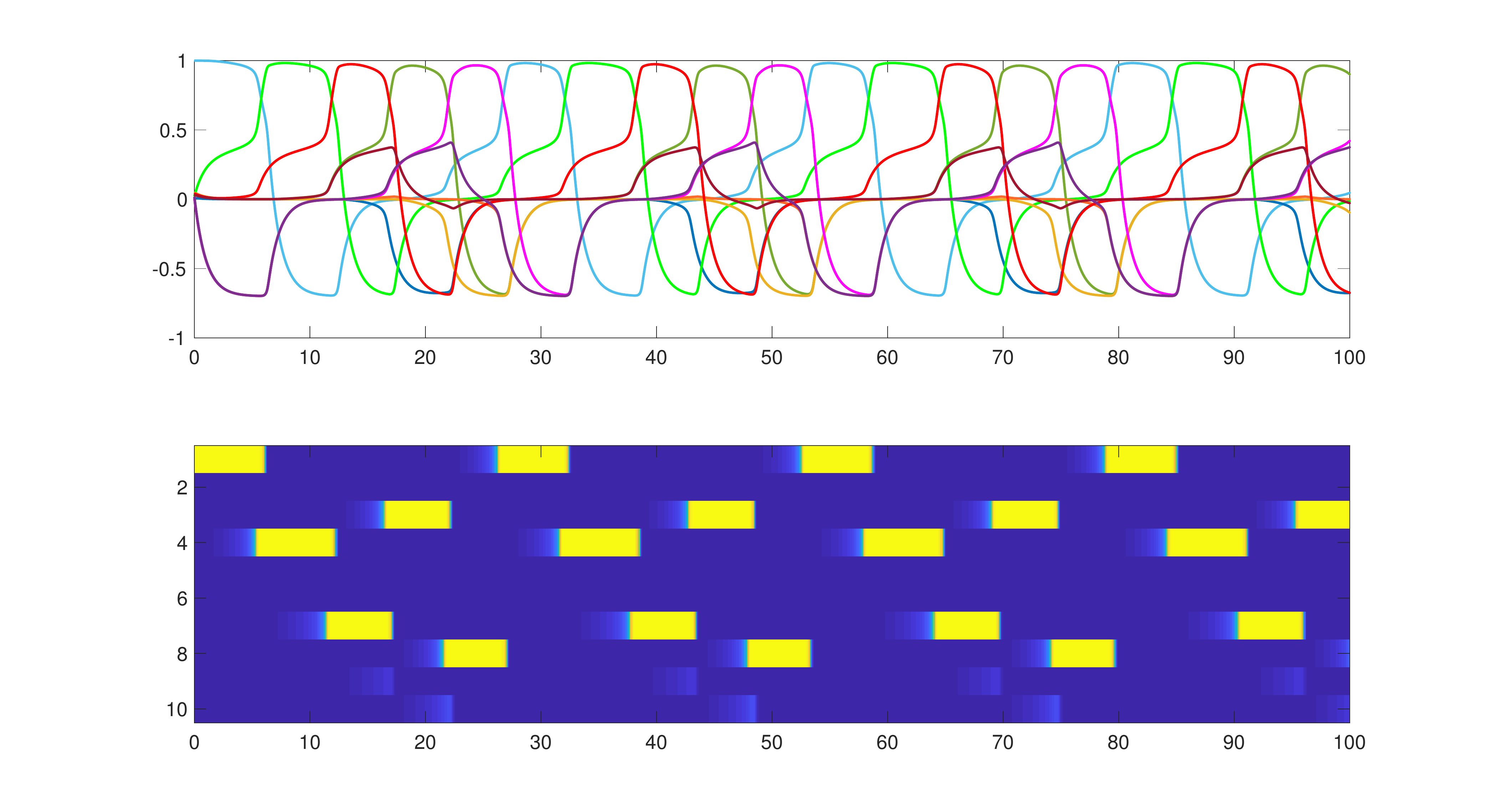}}
		\put(135,5){$t$}
		\put(0,80){$y_j$}
		\end{picture}
	\end{center}
	\caption{Trajectories of the CTRNN network for the directed graph shown in figure~\ref{fig:tennodegraph}, using the deterministic model~\eqref{eq:ctrnn_ode}. The top panel shows the $y_j$ coordinates, where the colours correspond to the node colours in figure~\ref{fig:tennodegraph}. In the lower panel, each horizontal row corresponds to one node (as labelled on the vertical axis), and the colour is blue when the corresponding $J_j$ coordinate is close to zero, and yellow when it is close to one. That is, the yellow segments indicate when each node is active. Parameters are as described in the text.
		\label{fig:tennodetraj_det}}
\end{figure*}

\begin{figure*}
	\setlength{\unitlength}{1mm}
	\begin{center}
		\begin{picture}(140,85)(0,5)
		\put(0,5){\includegraphics[trim= 3.5cm 0.5cm 3.2cm 0.5cm,clip=true,width=14cm]{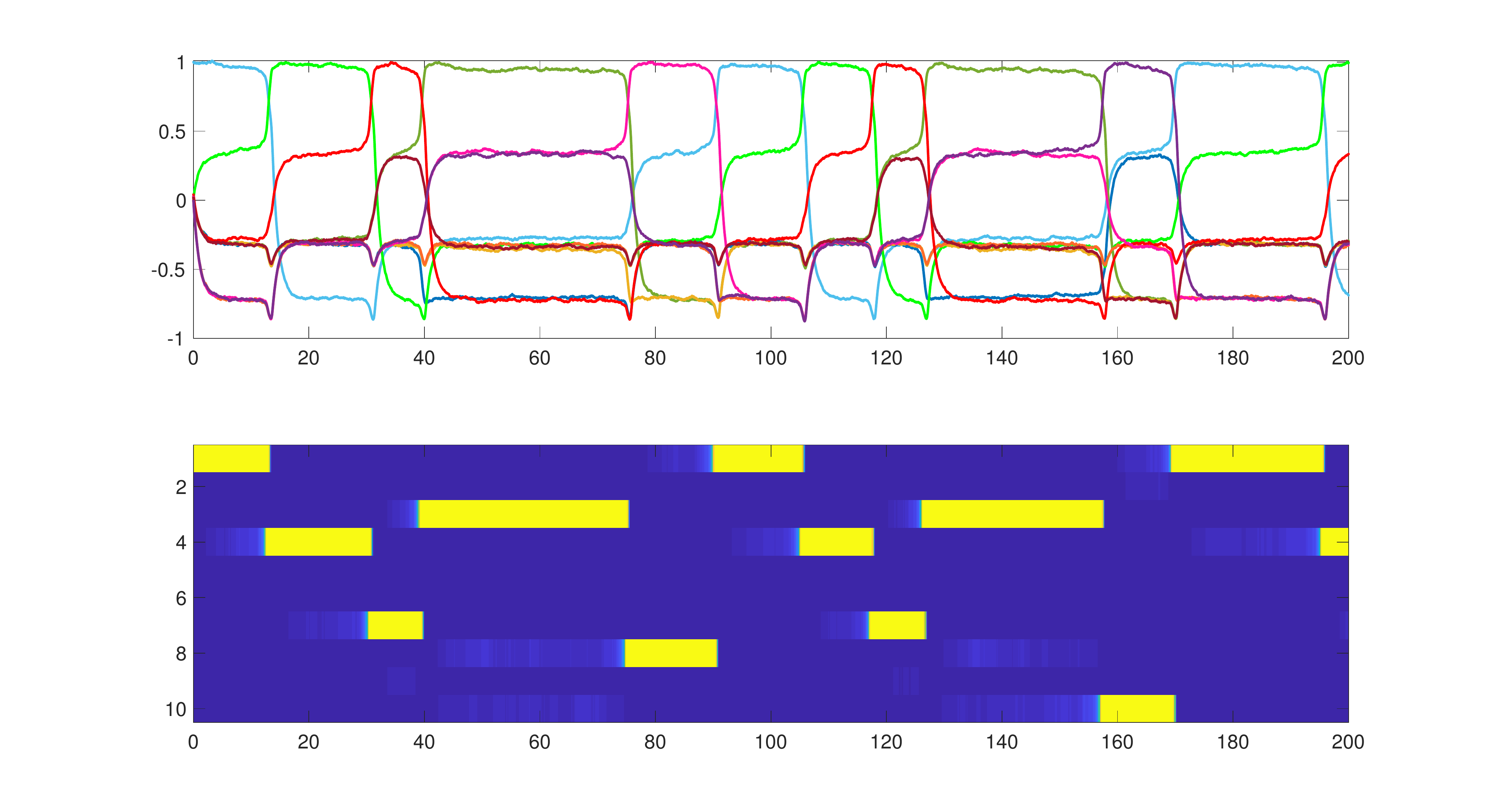}}
		\put(135,5){$t$}
		\put(0,80){$y_j$}
		\end{picture}
	\end{center}
	\caption{Trajectories of the CTRNN network for the directed graph shown in figure~\ref{fig:tennodegraph}, using the stochastic model~\eqref{eq:ctrnn_sde}. Lines and colours are described in figure~\ref{fig:tennodetraj_det}. Parameters are as described in the text, here $w_t=-0.3$.
		\label{fig:tennodetraj_noisy1}}
\end{figure*}

\begin{figure*}
	\setlength{\unitlength}{1mm}
	\begin{center}
		\begin{picture}(140,85)(0,5)
		\put(0,5){\includegraphics[trim= 3.5cm 0.5cm 3.2cm 0.5cm,clip=true,width=14cm]{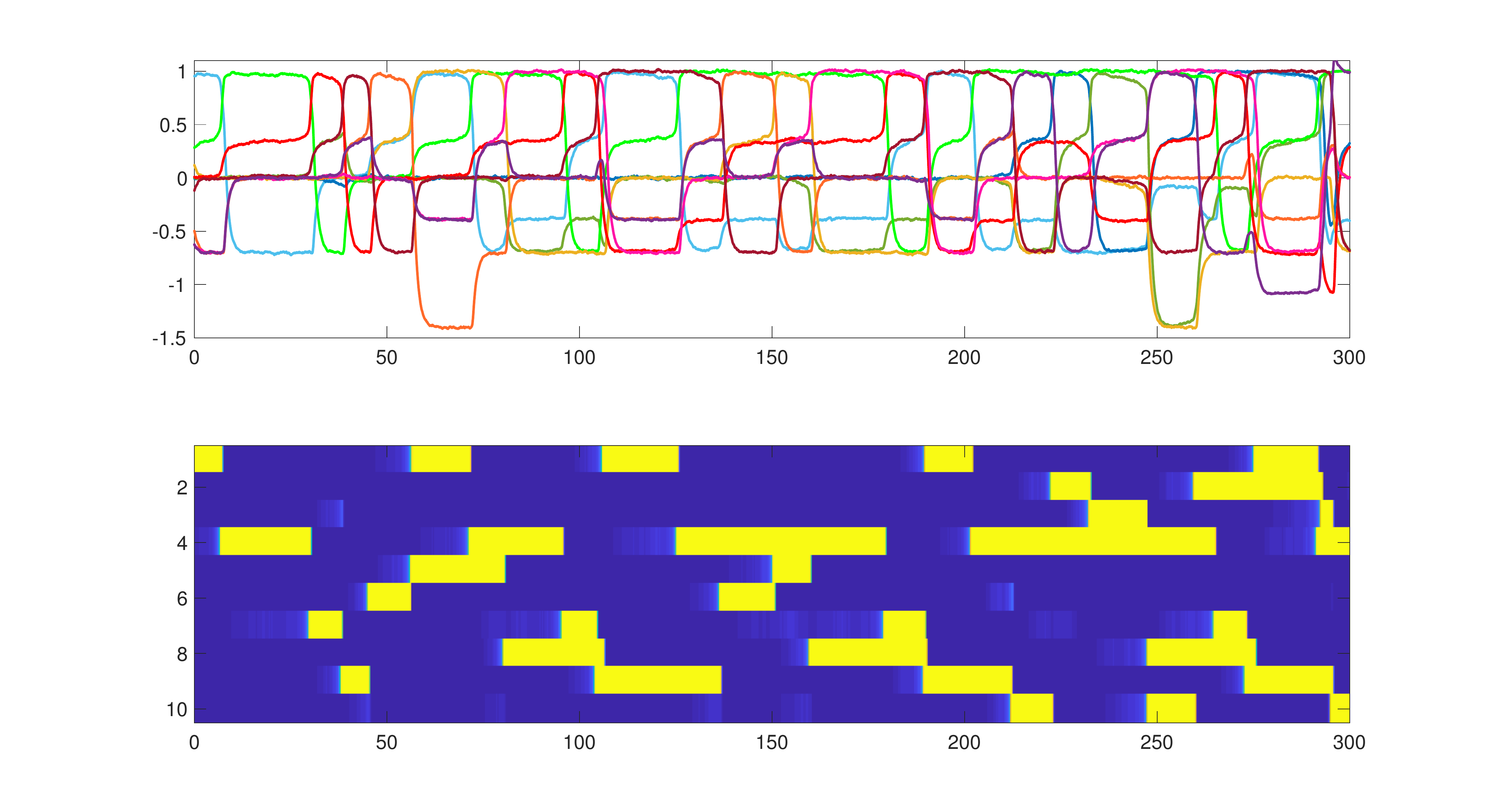}}
		\put(135,5){$t$}
		\put(0,80){$y_j$}
		\end{picture}
	\end{center}
	\caption{Trajectories of the CTRNN network for the directed graph shown in figure~\ref{fig:tennodegraph}, using the stochastic model~\eqref{eq:ctrnn_sde}. Lines and colours are described in figure~\ref{fig:tennodetraj_det}. Parameters are as described in the text, here $w_t=0$.
		\label{fig:tennodetraj_noisy2}}
\end{figure*}

The results of the simulations are shown in figures~\ref{fig:tennodetraj_det}, \ref{fig:tennodetraj_noisy1} and \ref{fig:tennodetraj_noisy2}. For the deterministic simulation (figure~\ref{fig:tennodetraj_det}), we see that the system has an attracting period orbit, which visits the nodes in the order $1\rightarrow 4 \rightarrow 7 \rightarrow 3 \rightarrow 8$. The entries of $w_p^{ij}$ were randomly generated as described above, and we do not give them all here for space reasons, but we note that in all cases in which a vertex in this cycle has two `choices' for which direction to leave (in the graph shown in figure~\ref{fig:tennodegraph}), the attracting periodic orbit chooses the more unstable direction. That is, if $i$ is a vertex in the above cycle, and if $i\rightarrow j$ and $i\rightarrow k$ are connections in the directed graph, with $i\rightarrow j$ being a part of the attracting periodic orbit, then $w_p^{ji}>w_p^{ki}$. Although we do not prove here that this will always be the case, it is intuitively what one might expect, that is, the connection from $i$ to $j$ is stronger than the connection from $i$ to $k$.

In the noisy simulation with $w_t=-0.3$ (figure~\ref{fig:tennodetraj_noisy1}),  the equilibria are not visited in a regular pattern, but random choices are made at each equilibria from which there is more than one direction in which to leave. See, for instance, the transition $3\rightarrow 8$ at $t\approx 75$, and the transition $3\rightarrow 10$ at $t\approx 155$. The length of time spent near each equilibria is also irregular; note for instance, the variable amount of time spent near $\xi_1$ and $\xi_3$. In this simulation, because the transverse parameter $w_t$ is sufficiently negative, only one node is active at any given time. 

By contrast, the simulation in figure~\ref{fig:tennodetraj_noisy2} has $w_t=0$. Here, the transitions again are made randomly, but without the suppression provided by the transverse parameter it is possible to have multiple cells active at once. For instance, at around $t=55$, both cells 1 and 5 become active at the same time; they were both leading cells from the previously active cell 6. As cells 1 switches off,  cell 4 becomes active, and as cell 5 switches off, cell 8 becomes active. The system continues to have two active cells around $t=250$, at which point a third cell also becomes active. If the trajectory were to run for longer then the number of active cells could decrease again, if an active cell suppresses more than one previously active cells.

The entire excitable network attractor for this level of noise is clearly more complicated than the design shown in figure~\ref{fig:tennodegraph}, in that additional equilibria (with more than one active cell) are accessible to those encoded and described in Theorem~\ref{thm:main}. An interesting extension of this work would be to understand which additional equilibria appear in a network attractor generated from a given directed graph in this manner: figure~\ref{fig:tennodetraj_noisy2} suggests that at least seven levels of cell activity are needed to uniquely describe the states that can appear when more than one cell becomes ``active''. In particular, when a cell is Trailing to more than one Active cell, the value of $y_j$ for that cell is even lower than $Y_T$ (compare the top panel of figure~\ref{fig:tennodetraj_noisy2} with the schematic in figure~\ref{fig:transition_schematic}).

\section{Discussion}
\label{sec:discuss}

The main theoretical result of this paper is Theorem~\ref{thm:main}, which states that it is possible to design the connection weight matrix of a CTRNN such that there exists a network attractor with a specific graph topology embedded within the phase space of the CTRNN. The graph topology is arbitrary except for minor restrictions: namely there should be no loops of order one or two, and no $\Delta$-cliques. Theorem~\ref{thm:main} assumes a piecewise affine activation function, but the examples in Section~\ref{sec:ex} suggest that the results generalise to CTRNN using any suitable smooth activation function. More generally, note that the coupled network is in some sense close to $N$ simultaneous saddle-node bifurcations. However, the units are not weakly coupled and indeed this is necessary to ensure that when one cell becomes active, the previous active cell is turned off.

Theorem~\ref{thm:main} proves the existence of an excitable network with threshold $\delta$ where not only the connection weights, but also $\epsilon$ and $\theta$ (properties of the activation function) may depend on $\delta$. We believe that a stronger result will be true, namely that $\delta$ can be chosen independent of properties of the activation function, and also that this can be made an almost complete realisation by appropriate choice of parameters.

\begin{conjecture}
\label{conj:main}
Assume the hypotheses on the directed graph $G$ with $N$ vertices as in Theorem~\ref{thm:main} hold. Assume that $\epsilon>0$ is small and $\theta>0$. Then there is a $\delta_c(\epsilon,\theta)>0$ such that for any $0<\delta<\delta_c$ there is an open set $\hat{W}_{\mathrm{ex}}\subset\R^{4}$. If the parameters $(w_s,w_m,w_t,w_p)\in \hat{W}_{\mathrm{ex}}$ then the dynamics of input-free equation~\eqref{eq:ctrnn_ode} with $N$ cells and piecewise affine activation function \eqref{eq:phip} and $w_{ij}$ defined by~\eqref{eq:w} contains an excitable network attractor with threshold $\delta$ that gives an almost complete realisation of the graph $G$.
\end{conjecture}

The construction in Theorem~\ref{thm:main} uses a comparatively sparse encoding of network states - each of the $N$ vertices in the network is associated with precisely one of the $N$ cells being in an active state. Indeed, the connection weights \eqref{eq:w} assign one of only four possible weights to each connection, depending on whether that cell can become active next, was active previously or neither. Other choices of weights will allow more dense encoding: and many more than $N$ excitable states within a network of $N$ cells. However the combinatorial properties of the dynamics seem to be much more difficult to determine and presumably additional connection weightings will be needed, not the just four values considered in Theorem~\ref{thm:main}.

Section~\ref{sec:ex} illustrates specific examples of simple excitable networks  for smooth activation function (\ref{eq:phi}) on varying parameters - this requires numerical continuation to understand dependence on parameters even for fairly low dimension. For this reason we expect that a proof of an analogous result to Theorem~\ref{thm:main} for the smooth activation function may be a lot harder. These examples also give some insight into bifurcations that create the excitable networks. 

In general there is no reason that the realisation constructed in Theorem~\ref{thm:main} is \emph{almost complete} (in the sense that almost all initial conditions in $B_{\delta}(\xi_k)$ evolve towards some $\xi_l$ with $a_{kl}=1$, analogous to \cite[Defn 2.6]{ashwin2020almost}). If it is not then other attractors may be reachable from the excitable network. Conjecture~\ref{conj:main} suggests that the realisation can be made almost complete for small enough $\delta$: we expect that for this we will require $w_t$ to be sufficiently negative. If $\delta$ is too large we cannot expect almost completeness: there are other stable equilibria (notably the origin) that can be reached with large perturbations, and the simulation in figure~\ref{fig:tennodetraj_noisy2} shows that other equilibria may be reachable from the network. It will be a challenge to strengthen our results to show that the excitable network is an almost complete realization. However, the examples studied in section~\ref{sec:ex} confirm that, at least for relatively simple graphs, this conjecture is reasonable. 

Our theoretical results are for networks with excitable connections. We expect much of the behaviour described here is present in the spontaneous case, if the coupling weights are chosen such that equilibria are replaced with bottlenecks. In the absence of noise we expect to see a deterministic switching between slow moving dynamics within bottlenecks. This dynamical behaviour is very reminiscent of the stable heteroclinic channels described, for example, in \cite{afraimovich2004heteroclinic,rabinovich2020sequential}. However, stable heteroclinic switching models require structure in the form of multiplicative coupling or symmetries that are not present in CTRNN or related Wilson-Cowan neural models \cite{wilson1972excitatory} (see \cite{chow2019before} for a recent review of related neural models). Other models showing sequential excitation include \cite{chow2019before}: this relies on a fast-slow decomposition to understand various different modes of sequential activation in a neural model of rhythm generation. It will be an interesting challenge to properly describe possible output dynamics of our model in the case of bottlenecks.

We remark that asymmetry of connection weights is vital for constructing a realisation as an excitable networks - indeed, the lack of two-cycles precludes $a_{jk}=a_{kj}=1$. While this may be intuituively obvious, it was not so obvious that we also need to exclude one-cycles and $\Delta$-cliques in the graph to make robust realisations. 

Finally, although we do not consider specific natural or machine learning applications of CTRNN here, the structures found here may give insights that give improved training for CTRNN. In particular, it seems plausible that CTRNN may use excitable networks to achieve specific input-output tasks (especially those requiring internal states). For example, recent work \cite{ceni2020interpreting} demonstrates that echo state networks can create excitable networks in their phase space to encode input-dependent behaviour. It is also likely there are novel optimal training strategies that take advantage of excitable networks, for example, choosing connection weights that are distributed close to one of the four values we use.

\begin{acknowledgements}
We gratefully acknowledge support from the Marsden Fund Council from New Zealand Government funding, managed by Royal Society Te Ap\={a}rangi. PA acknowledges funding from EPSRC as part of the Centre for Predictive Modelling in Healthcare grant EP/N014391/1. CMP is grateful for additional support from the London Mathematical Laboratory.
\end{acknowledgements}

\appendix

\section{Definition of excitable network}

%x
\label{app:networks}

This appendix (which extends ideas in \cite{Ashwin2016a}) gives formal definitions for excitable networks considered in this paper. We say a system has an \emph{excitable connection for amplitude} $\delta>0$ from one equilibrium $\xi_i$ to another $\xi_j$ if 
\[
B_{\delta}(\xi_i)\cap W^s(\xi_j)\neq \emptyset,
\]
(where $B_{\delta}(\xi)$ is the open ball of radius $\delta$ centred at $\xi$) and this connection has \emph{threshold} $\delta_{th}$ if
\[
\delta_{th} =\inf\{\delta>0 : B_{\delta}(\xi_i)\cap W^s(\xi_j)\neq \emptyset\}.
\]
We define the \emph{excitable network of amplitude} $\delta>0$ between the equilibria $E=\{\xi_i\}$ to be the set
\[
\Sigma_E=\bigcup_{i,j=1}^N\{\Phi_t(x): x\in B_{\delta}(\xi_i), t>0 \} \cap W^s(\xi_j).
\]

We say the excitable network $\Sigma_E$ for amplitude $\delta$ \emph{realises a graph} $G$ if each vertex $v_i$ in $G$ corresponds to an equilibrium $\xi_i$ in $\Sigma_E$ and there is an edge in $G$ from $v_i$ to $v_j$ if only if there is a connection in $\Sigma_E$ for amplitude $\delta$ from $\xi_i$ to $\xi_j$. 
%We say $\Sigma_E$ \emph{properly realises} $G$ when there is an excitable connection from $\xi_i$ to $\xi_j$ in $E$ if and only if there is a corresponding edge in $G$.

\section{Proof of Theorem~\ref{thm:main} }
\label{app:proof}

As in equations~\eqref{eq:paramchoicea} and~\eqref{eq:paramchoiceb}, for any $\delta<\frac{1}{2}$ we choose 
\begin{equation}
\epsilon=\dfrac{\delta}{8},~\theta=\frac{1}{2},~w_s=1,~w_t=0,
\end{equation}
and then $w_p$ and $w_m$ are given by
\begin{equation}
w_p=\theta-\dfrac{\delta}{2},\quad w_m=-(w_s-\theta)-\dfrac{\delta}{2}.
\end{equation}
We define equilibria $\xi_k$ as in Section~\ref{sec:thm}, where we write the $j$th component of the equilibrium as $[\xi_k]_j$:
\begin{equation}
\label{eq:xikj}
[\xi_k]_j=\begin{cases}
Y_A & \mbox{ if }j=k,\\
Y_L & \mbox{ if }a_{kj}=1,\\
Y_T & \mbox{ if }a_{jk}=1,\\
Y_D & \mbox{ if }a_{kj}=0\mbox{ and }a_{jk}=0,
\end{cases}
\end{equation}
where 
\begin{align*}
Y_A&:=w_s=1,\ Y_L:=w_p=\theta-\frac{\delta}{2},\\ 
Y_T&:=w_m=-(w_s-\theta)-\frac{\delta}{2},\ Y_D:=w_t=0,
\end{align*}
are the values of the Active, Leading, Trailing and Disconnected components, respectively. 
Note that
\begin{equation}
\label{eq:ineq}
Y_T<Y_D<Y_L<\theta-\epsilon/2<\theta+\epsilon/2<Y_A.
\end{equation}
As mentioned before, the hypotheses of Theorem~\ref{thm:main} imply that this labelling is well defined and Figure~\ref{fig:transition_schematic} shows how a transition from cell 1 active to 
cell 2 active will occur in a general network.  It is simple to check that~\eqref{eq:xikj} is an equilbrium solution of \eqref{eq:ctrnn_ode} with activation function~\eqref{eq:phip}. Moreover, $\xi_k$ is linearly stable with $n$ eigenvalues $-1$.
We define  $[\mathcal{J}_k]_j:=\phi_P([\xi_k]_j)$ then note that \eqref{eq:ineq} implies that
$$
[\mathcal{J}_k]_j:=\delta_{kj},
$$
in terms of the Kronecker $\delta_{kj}$.

\begin{figure}

	\setlength{\unitlength}{0.75mm}
	\begin{center}
		\begin{picture}(105,85)(0,0)
	\put(0,0){\includegraphics[width=7.5cm]{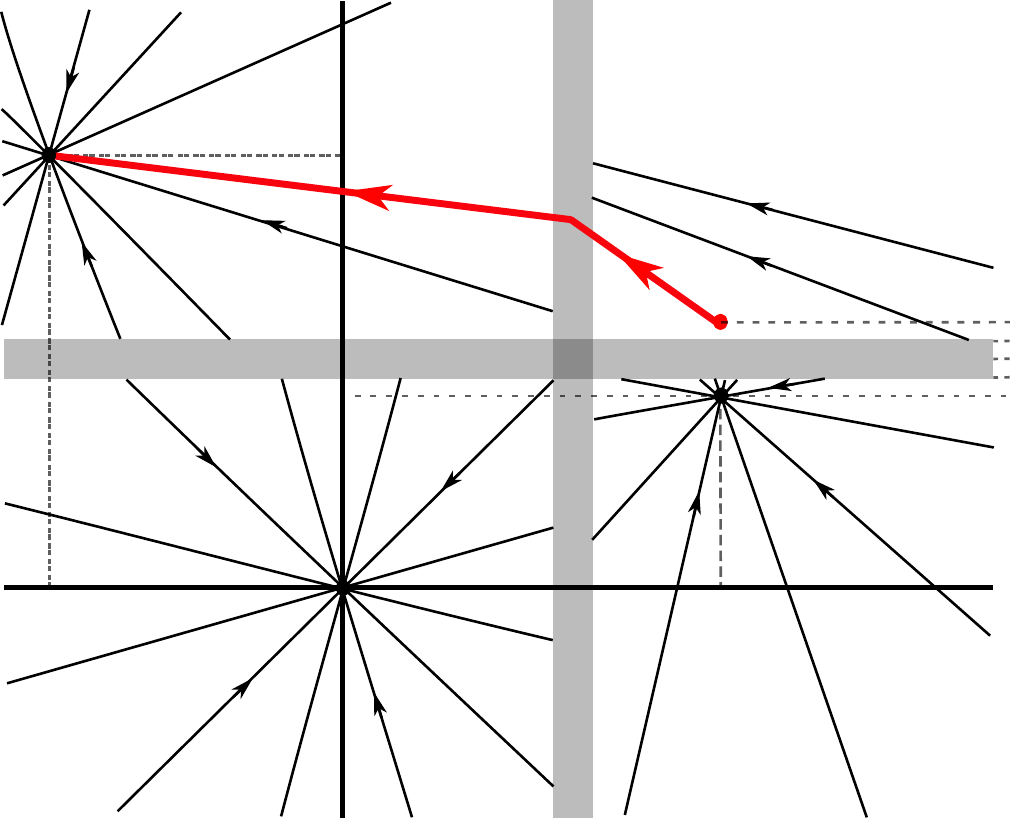}}
	
	\put(98,20){$y_k$}
	\put(30,82){$y_l$}
\put(80,70){$\mathcal{R}_{kl}$}  
\put(40,70){$\mathcal{R}_l$}
\put(55,83){$\mathcal{S}_{kl}$}
\put(88,30){$\mathcal{R}_k$}

\put(64,37.5){$\xi_k$}
\put(70,52){$\zeta_{k,l}$}
\put(3.5,73){$\xi_l$}
	
	\put(71,20){$w_s$}
	\put(29,41){$w_p$}
	
	\put(29.5,68){$w_s$}
	\put(3,20){$w_m$}
	
	\put(31,20){$0$}
	
	\put(100.5,41){\tiny $\theta-4\epsilon$}	
	\put(100.5,43){\tiny $\theta-2\epsilon$}	
	\put(100.5,45){\tiny $\theta$}
	\put(100.5,47){\tiny $\theta+2\epsilon$}	
	\put(100.5,49){\tiny $\theta+4\epsilon$}

%	\put(0,0){x}
%	\put(100,85){x}
	
	\end{picture}
	\end{center}
	
	\caption{Illustration of a connection (shown in red) with amplitude $\delta>0$ from $\xi_k$ to $\xi_l$, projected in the coordinates $y_k$ and $y_l$. The points $\xi_k$, $\xi_l$ and $0$ are linear sinks. The point $\zeta_{k,l}$ is within $\delta$ of $\xi_l$ and limits to $\xi_l$ in forwards time. The grey areas are regions of width $4\epsilon$ centred around $y_{k,l}=\theta$ where $\phi_P(y_k)$ and $\phi_P(y_l)$ are non-constant: these contain saddle and other equilibria (not shown).}
	\label{fig:connection_phipwl}
\end{figure}
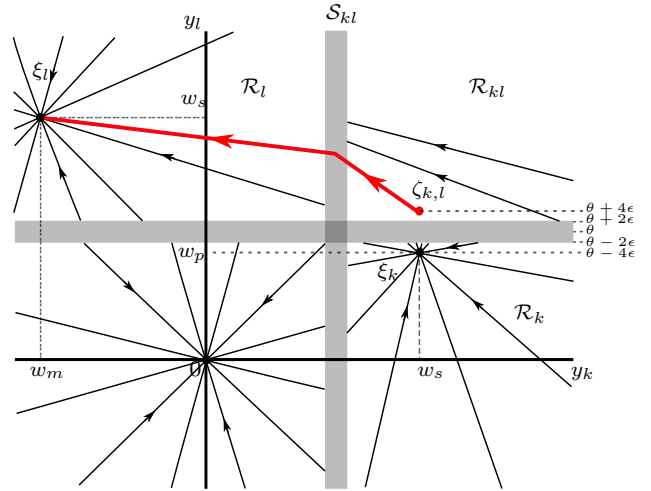

We consider two cases. Case 1 is where $a_{kl}=1$ and we expect to see a connection from $\xi_k$ to $\xi_l$. Case 2 is $a_{kl}=0$ and we don't expect a connection.

\subsection{Case 1}

Suppose $a_{kl}=1$ (recall that the lack of $1$-cycles means that $k\neq l$), and then we define two regions of phase space, which we label $\mathcal{R}_{l}$ and $\mathcal{R}_{kl}$, as follows:
 \[
 \mathcal{R}_{l} =\left\{ \vec{y}\ |\ y_j>\theta+\frac{\delta}{4},\ j=l;\ y_j<\theta-\frac{\delta}{4}\ \text{otherwise}\right\},
\]
 \[
 \mathcal{R}_{kl} =\left\{ \vec{y}\ |\ y_j>\theta+\frac{\delta}{4},\ j=l,k;\ y_j<\theta-\frac{\delta}{4}\ \text{otherwise}\right\}.
\]
The regions and the dynamics within them are shown schematically in figure~\ref{fig:connection_phipwl}.

Within $\mathcal{R}_{kl}$, if $a_{kl}=1$ then the dynamics are governed by the equations
\[
\dot{y}_j=f_j^{(b)}:=-y_j+\phi^{(b)}_j,
\]
where
\begin{equation}
\label{eq:xdotcasesjkl}
\phi_j^{(b)}:=\begin{cases}
%\underline{\mbox{ If }j\in\{k,l\}:}&\\
w_s+w_m & \mbox{ if }j=k, \text{(AT)}\\
w_s+w_p & \mbox{ if }j=l,  \text{(LA)}\\
%\underline{\mbox{ If }j\not\in\{k,l\}:}&\\
w_p & \mbox{ if } (a_{kj},a_{jk},a_{lj},a_{jl})=(1,0,0,0),  \text{(LD)}\\
w_m & \mbox{ if } (a_{kj},a_{jk},a_{lj},a_{jl})=(0,1,0,0),  \text{(TD)}\\
w_p & \mbox{ if } (a_{kj},a_{jk},a_{lj},a_{jl})=(0,0,1,0),  \text{(DL)}\\
w_m & \mbox{ if } (a_{kj},a_{jk},a_{lj},a_{jl})=(0,0,0,1),  \text{(DT)}\\
w_m + w_p & \mbox{ if } (a_{kj},a_{jk},a_{lj},a_{jl})=(0,1,1,0),  \text{(TL)}\\
0  & \mbox{ if } (a_{kj},a_{jk},a_{lj},a_{jl})=(0,0,0,0),  \text{(DD)}
\end{cases}
\end{equation}
where the type of coordinate (see figure~\ref{fig:transition_schematic}) is given in parentheses on each line. 
Within $\mathcal{R}_{l}$, the dynamics are governed by the equations
\begin{equation}
\label{eq:xdotcases_jl}
\dot{y}_j=f_j^{(c)}:=-y_j+\phi^{(c)}_j,
\end{equation}
where
\begin{equation}
\phi_j^{(c)}:=\begin{cases}
w_m & \mbox{ if }j=k, \text{(AT)} \\
w_s & \mbox{ if }j=l,  \text{(LA)}\\
w_p & \mbox{ if }a_{lj}=1 \mbox{ and }j\neq l,  \text{(DL/TL)}\\
w_m & \mbox{ if }a_{jl}=1, \text{(DT)}\\
0  & \mbox{ if }a_{lj}=0\mbox{ and }a_{jl}=0, \text{(LD/TD/DD)}
\end{cases}
\end{equation}

%\pa{Edited next sentence.. please check}
The equilbrium $\xi_l$ lies in the interior of the region $\mathcal{R}_l$, whereas $\xi_k$ lies in the interior of the region $\mathcal{R}_k$. We show 
 that there is a connection of amplitude $\delta$ from $\xi_k$ to any $\xi_l$ with $a_{kl}=1$, by considering a trajectory starting at 
$$
\zeta_{k,l}=\xi_k+\delta e_{l},
$$
where $e_l$ is a unit vector in the $l$-direction: clearly $|\zeta_{k,l}-\xi_k|=\delta$ (see figure~\ref{fig:connection_phipwl}). Note that
\begin{equation}
\label{eq:xitilde}
[\zeta_{k,l}]_j=\begin{cases}
Y_A & \mbox{ if }j=k,\\
Y_L+\delta & \mbox{ if } j=l,\\% \mbox{ (and so }a_{kl}=1\mbox{)}\\
Y_L & \mbox{ if }a_{kj}=1 \mbox{ and }j\neq l,\\
Y_T & \mbox{ if }a_{jk}=1,\\
Y_D & \mbox{ if }a_{kj}=0\mbox{ and }a_{jk}=0,
\end{cases}
\end{equation}
where $Y_A$ etc are given in equation~\eqref{eq:yaltd}, and $Y_L+\delta=\theta+\delta/2>\theta$.
We define $[\mathcal{J}_{k,l}]_j:=\phi_P([\zeta_{k,l}]_j)$, and then 
$$
[\mathcal{J}_{k,l}]_j=\delta_{jk}+\delta_{jl}.
$$

Now, we note that $\zeta_{k,l}\in\mathcal{R}_{kl}$. We next show that a trajectory with initial condition at $\zeta_{k,l}$ will asymptotically approach $\xi_l$. We show first that the trajectory enters $\mathcal{R}_l$ in finite time. Then, since in $\mathcal{R}_l$, the flow is linear, with stable equilibrium $\xi_l$, all trajectories in $\mathcal{R}_l$ eventually approach $\xi_l$.

Define $\mathcal{S}_{kl}$ to be the region between $\mathcal{R}_l$ and $\mathcal{R}_{kl}$, namely,
 \[
 \mathcal{S}_{kl} =\left\{ \vec{y}\ |\ \theta-\frac{\delta}{4}\leq y_k \leq \theta+\frac{\delta}{4}; y_l>\theta+\frac{\delta}{4} ;\ y_j<\theta-\frac{\delta}{4}\ \text{ for } j\neq k,l\right\}.
\]

First consider the dynamics of all $y_j$, with $j\neq k,l$. This means that $y_j<\theta-\delta/4$ for all points in $\mathcal{T}_{kl}\equiv \mathcal{R}_{kl} \cup \mathcal{S}_{kl} \cup \mathcal{R}_l$. While the trajectory remains in $\mathcal{T}_{kl}$, it can be shown that $\dot{y}_j$ is negative along the line with $y_j=\theta-\delta/4$. Hence, none of the $y_j$ will leave $\mathcal{T}_{kl}$.

Within $\mathcal{S}_{kl}$, the dynamics of $y_k$ and $y_l$ are governed by
\begin{align}
\dot{y}_k=&-y_k+w_m+\phi_P(y_k), \\
\dot{y}_l=&-y_l+w_s+w_p\phi_P(y_k). 
\end{align}

Consider the equation for $\dot{y}_k$ in $\mathcal{R}_{kl}$, $\mathcal{S}_{kl}$ and $\mathcal{R}_l$, namely:
\[
\dot{y}_k=\begin{cases} 
-y_k+w_s+w_m & \text{if}\ \vec{y}\in  \mathcal{R}_{kl}, \\ 
-y_k+\phi_P(y_k)+w_m & \text{if}\ \vec{y}\in  \mathcal{S}_{kl},\\ 
-y_k+w_m & \text{if}\ \vec{y}\in  \mathcal{R}_{l}. 
\end{cases}
\]
Recall that $w_m=-(w_s-\theta)-\frac{\delta}{2}$, and $0\leq \phi_P(y_k) \leq 1$. We can use these bounds to show that
\[
\dot{y}_k \leq \begin{cases} 
-\dfrac{3\delta}{4} & \text{if}\ \vec{y}\in  \mathcal{R}_{kl}, \\ 
-\dfrac{\delta}{4} & \text{if}\ \vec{y}\in  \mathcal{S}_{kl}, 
\end{cases}
\]
Furthermore, if $\vec{y}\in  \mathcal{R}_{l}$ and $y_k>0$, then  $\dot{y}_k<-w_m<0$. In particular, we note that for $\vec{y}\in\mathcal{T}_{kl}$, with $y_k>0$, $\dot{y}_k$ is negative and bounded below zero. 

Now, consider the equation for $\dot{y}_l$ in $\mathcal{R}_{kl}$, $\mathcal{S}_{kl}$ and $\mathcal{R}_l$, namely:
\[
\dot{y}_l=\begin{cases} 
-y_l+w_s+w_p & \text{if}\ \vec{y}\in  \mathcal{R}_{kl}, \\ 
-y_l+w_s+w_p\phi_P(y_k) & \text{if}\ \vec{y}\in  \mathcal{S}_{kl}, \\ 
-y_l+w_s & \text{if}\ \vec{y}\in  \mathcal{R}_{l}, 
\end{cases}
\]
We use this to compute $\dot{y}_l$ along the lower boundary of the three regions, $\mathcal{R}_{kl}$, $\mathcal{S}_{kl}$ and $\mathcal{R}_l$, that is, the line $y_l=\theta+\frac{\delta}{4}$, and we find
\[
\dot{y}_l=\begin{cases} 
w_s-\frac{3\delta}{4}& \text{if}\ \vec{y}\in  \mathcal{R}_{kl}, \\ 
w_s-\theta-\frac{\delta}{4}+\left(\theta-\frac{\delta}{2} \right)\phi_P(y_k) >w_s-\frac{3\delta}{4}& \text{if}\ \vec{y}\in  \mathcal{S}_{kl}, \\ 
w_s-\theta-\frac{\delta}{4} & \text{if}\ \vec{y}\in  \mathcal{R}_{l}. 
\end{cases}
\]
Since $w_s>\frac{3\delta}{4}$ and $w_s>\theta+\frac{\delta}{4}$, we see that $\dot{y}_l>0$ in all three cases.

Combining our knowledge of $\dot{y}_l$ and $\dot{y}_k$ tells us that a trajectory which starts in $\mathcal{R}_{kl}$, or more specifically, a trajectory starting in a small neighbourhood of $\zeta_{k,l}$ will have monotonic decreasing $y_k$ component until (at least) $y_k=0$. Furthermore, the $y_l$ component cannot decrease below $y_l=\frac{1}{2}+\frac{\delta}{4}$. Thus the trajectory will move through $\mathcal{S}_{kl}$ and into $\mathcal{R}_{l}$ in a bounded time. 

Within $\mathcal{R}_{l}$, $\xi_l$ is a linearly stable fixed point.  In summary, we have demonstrated that if $a_{kl}=1$ then there is a connection from a $\delta$-neighbourhood of $\xi_k$ to $\xi_l$. Moreover, as the equilibria are linearly stable and having a connection is an open condition, the realisation will persist for an open set of parameters.

\subsection{Absence of excitable connections for edges absent from $G$}
\label{sec:proofproper}

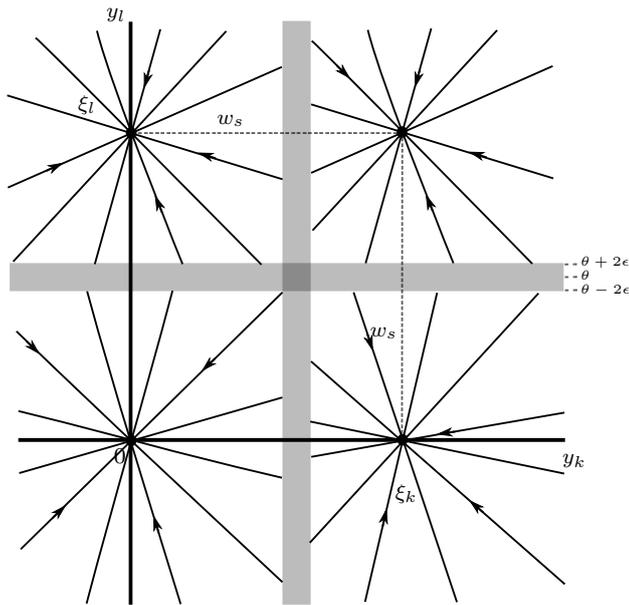
\begin{figure}

	\setlength{\unitlength}{0.9375mm}
	\begin{center}
		\begin{picture}(105,85)(0,0)
	\put(0,0){\includegraphics[width=7.5cm]{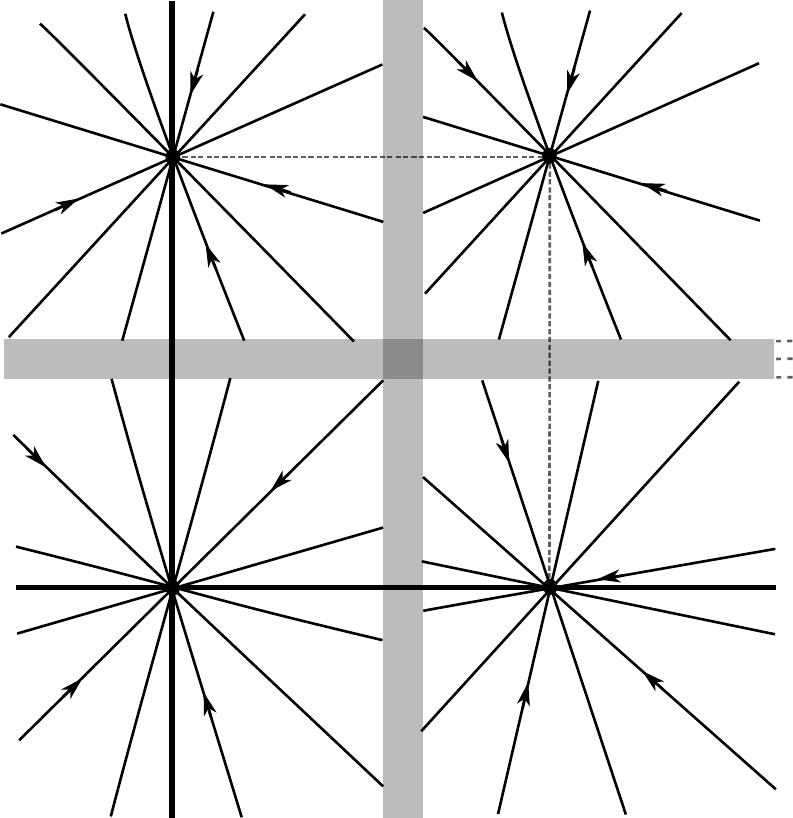}}
	
	\put(78,20){$y_k$}
	\put(14,83){$y_l$}

	\put(54.5,15){$\xi_k$}
	\put(10,70){$\xi_l$}
	
	\put(51,38){$w_s$}
	
	\put(29.5,68){$w_s$}
	
	\put(15,20){$0$}
	
	\put(80.5,44){\tiny $\theta-2\epsilon$}	
	\put(80.5,46){\tiny $\theta$}
	\put(80.5,48){\tiny $\theta+2\epsilon$}	
	
	\end{picture}
	\end{center}
	
	\caption{Schematic diagram showing the dynamics in the $y_k$-$y_l$ plane when $a_{kl}=a_{lk}=0$. Equilibria are shown with dots, and all equilibria are linearly stable. }
	\label{fig:noconnection_phipwl}
\end{figure}

Now suppose that $a_{lk}=a_{kl}=0$. Then the dynamics for $x_k$ and $x_l$ is shown schematically in figure~\ref{fig:noconnection_phipwl}.  Equilibria are shown with dots, and all equilibria shown in this figure are linearly stable.
Note that the equilibrium $\xi_k$ has $y_k=Y_A=1$, and $y_l=Y_D=0$. It is clear that there are no small perturbations which allow for a connection between $\xi_k$ and $\xi_l$. 

For $\theta=1/2$, if $\delta<1/4$ then for any $k$ and $j\neq k$ such that $a_{kj}=0$, all trajectories starting in
$$
\zeta_{k,l}=\xi_k+a e_k+b e_l,
$$
with $a^2+b^2\leq \delta^2$. In particular, perturbations of the form (\ref{eq:xitilde}) will return to $\xi_k$.  This is because this set is remains in the region of validity for the equivalent of (\ref{eq:xdotcases_jl}) for which the only attractor is $\xi_k$. In the case $a_{kl}=0$ and $a_{lk}=1$ a similar argument holds as the phase portrait corresponds to figure~\ref{fig:connection_phipwl} reflected in the diagonal.

This shows that no perturbations of amplitude $\delta$ within the $(x_k,x_l)$ plane that give a connection from $\xi_k$ to $\xi_l$ of amplitude $\delta$. However, we cannot rule out the existence of connections from other locations in $B_{\delta}(\xi_k)$ to $\xi_l$. This would be needed to prove that the realisation is almost complete.

 \bibliographystyle{spmpsci} 
 \bibliography{ctrnn_refs}

\end{document}